\documentclass[11pt]{article}

\usepackage{graphics,graphicx}

\usepackage{color}

\usepackage{latexsym,multicol,amssymb}
\usepackage{amsmath, amsbsy}
\usepackage{amsopn, amstext}

\newcommand{\bean}{\begin{eqnarray*}}
\newcommand{\eean}{\end{eqnarray*}}
\newcommand{\bea}{\begin{eqnarray}}
\newcommand{\eea}{\end{eqnarray}}
\newcommand{\nn}{\nonumber}
\newcommand{\EQ}{\begin{equation}\begin{array}{llllllllll}}
\newcommand{\EE}{\end{array}\end{equation}}

\newtheorem{theorem}{Theorem}
\newtheorem{remark}{Remark}[section]
\newtheorem{lemma}{Lemma}
\newtheorem{example}{Example}

\newtheorem{corollary}{Corollary}
\newtheorem{proposition}{Proposition}

\newcommand{\refB}{(\ref{state-eq})-(\ref{cons-c})}

\newcommand{\refBN}{(\ref{state-N})-(\ref{bound2})}
\newcommand{\NK}{^{Nk}}
\newcommand{\bx}{\bar x}
\newcommand{\bu}{\bar u}

\newcommand{\R}{\Re}
\newcommand{\bb}{{\boldsymbol b}}
\newcommand{\bd}{{\boldsymbol d}}
\newcommand{\BN}{${\rm B}^{\rm N}$}

\def\Fr{\displaystyle \frac}
\newcommand{\MT}{\left[\begin{array}{cccccccccccccc}}
\newcommand{\EM}{\end{array}\right]}
\def\dsum{\displaystyle \sum}
\newenvironment{proof}{{\it Proof. }}{\hfill $\Box$}

\begin{document}

\title{ {\bf On the Rate of Convergence for the Pseudospectral Optimal Control of Feedback Linearizable Systems}
\thanks{The research was supported in part by AFOSR and AFRL} }

\author{Wei Kang\\
Department of Applied Mathematics\\
Naval Postgraduate School\\
Monterey, CA 93943\\
 {\tt\small wkang@nps.edu}
}
\date {}
\maketitle

\begin{abstract}
Over the last decade, pseudospectral (PS) computational methods for nonlinear constrained optimal control have been applied to many industrial-strength problems, notably the recent zero-propellant-maneuvering of the International Space Station performed by NASA.
In this paper, we prove a theorem on the rate of convergence for the optimal cost computed using PS methods. It is a first proved convergence rate in the literature of PS optimal control. In addition to the high-order convergence rate, two theorems are proved for the existence and convergence of the approximate solutions. This paper contains several essential differences from existing papers on PS optimal control as well as some other direct computational methods. The proofs do not use necessary conditions of optimal control. Furthermore, we do not make coercivity type of assumptions. As a result, the theory does not require the local uniqueness of optimal solutions. In addition, a restrictive assumption on the cluster points of discrete solutions made in existing convergence theorems are removed.
\end{abstract}

\section{Introduction}

Despite the fact that optimal control is one of the oldest problems in the history of control theory, practical tools of solving nonlinear optimal control problems are limited. Preferably, a feedback control law is derived from a solution to the famously difficult Hamilton-Jacobi-Bellman (HJB) equation. However, analytic solutions of this partial differential equation can rarely be found for systems with nonlinear dynamics. Numerical approximation of such solutions suffers from the well-known curse of dimensionality and it is still an open problem for systems with moderately high dimension. A practical alternative is to compute one optimal trajectory at a time so that the difficulty of solving HJB equations is  circumvented. Then this open-loop optimal control can be combined with an inner-loop tracking controller; or it can be utilized as a core instrument in a real-time feedback control architecture such as a moving horizon feedback. A critical challenge in this approach is to develop reliable and efficient computational methods that generate the required optimal trajectories. In this paper, we focus on some fundamental issues of pseudospectral computational optimal control methods. 

As a result of significant progress in large-scale computational algorithms and nonlinear programming, the so-called direct computational methods have become popular for solving
nonlinear optimal control problems \cite{betts:book, betts:survey, polak:book},
particularly in aerospace applications \cite{paris:aas06,
riehl:aas06}. In simple terms, in a direct method, the
continuous-time problem of optimal control is discretized, and the
resulting discretized optimization problem is solved by nonlinear
programming algorithms. Over the last decade, pseudospectral (PS) methods have emerged as a popular direct methods for optimal control. They have been applied to many industrial-strength problems, notably the recent attitude maneuvers of the International Space Station performed by NASA. By following an attitude trajectory developed using PS optimal control, the International Space Station (ISS) was maneuvered 180 degrees on March 3, 2007, by using the gyroscopes equipped on the ISS without propellant consumption. This single maneuver have saved NASA about one million dollars' worth of fuel \cite{kangbedrossian}. The Legendre PS optimal control method has already been developed into software named DIDO, a MATLAB based package commercially available \cite{ross}. In addition, the next generation of the OTIS software package \cite{otis} will have the Legendre PS method as a problem solving option. 

PS methods have been widely applied in scientific computation for models governed by partial differential equations. The method is well known for being very efficient in approximating solutions of differential equations. However, despite its success and several decades of development, the intersection between PS methods and nonlinear optimal control becomes an active research area only after the mid-1990's (\cite{elnagar2,fahroo}). As yet, many fundamental theoretical issues are still widely open. For the last decade, active research has been carried out in the effort of developing a theoretical foundation for PS optimal control methods. Among the research focuses, there are three fundamental issues, namely the state and costate approximation, the existence and convergence of approximate solutions, and the convergence rate. The general importance of these issues is not limited to PS methods. They are essential to other computational methods suchlike those based on Euler \cite{DH2000} and Runge-Kutta \cite{H2000} discretization. Similar to other direct computational optimal control methods, PS method are based upon the Karush-Kuhn-Tucker (KKT) conditions rather than the Pontryagin's Minimum Principle (MPM). In \cite{fahroo} and \cite{gong2}, a covector mapping was derived between the costate from KKT condition and the costate from PMP. The covector mapping facilitates a verification and validation of the computed solution. For the problem of convergence, some theorems were published in \cite{gong}; and then the results were generalized in \cite{kang2} to problems with non-smooth control. 

Among the three fundamental issues mentioned above, the most belated activity of research is on the rate of convergence. In fact, there have been no results published on the convergence rate for PS optimal control methods. Although some results on the issue of convergence were proved in \cite{gong} and \cite{kang2}, a main drawback of these results is the strong assumption in which the derivatives of the discrete approximate solutions are required to converge uniformly. In this paper, we prove a rate of convergence for the approximate optimal cost computed using PS methods. Then, we prove theorems on existence and convergence without the restrictive assumption made in \cite{gong} and \cite{kang2}. In addition to the high-order convergence rate addressed in Section \ref{rate}, which is the first proved convergence rate in the literature of PS optimal control, this paper contains several essential differences from existing papers on PS optimal control as well as some other direct computational methods. First of all, the proof is not based on necessary conditions of optimal control. Furthermore, we do not make coercivity type of assumptions. As a result, the theory does not require the local uniqueness of optimal solutions. Therefore, it is applicable to problems with multiple optimal solutions. Secondly, the proof is not build on the bases of consistent approximation theory \cite{polak:book}. Thus, we can remove the assumption in \cite{gong} and \cite{kang2} on the existence of cluster points for the derivatives of discrete solutions. The key that makes these differences possible is that we introduce a set of sophisticated regularization conditions in the discretization so that the computational algorithm has a greater control of the boundedness of the approximate solutions and their derivatives. Different from the existing results in the literature of direct methods for optimal control, the desired boundedness is achieved not by making assumptions on the original system, but by implementing specially designed search region for the discrete problem of nonlinear programming. This new boundary of search region automatically excludes possible bad solutions that are numerically unstable. 

The paper is organized as follows. In Section \ref{formulation}, the formulations of the optimal control problem and its PS discretization are introduced. In \ref{rate}, we prove two theorems on the rate of convergence. In Section \ref{sec-convergence}, two theorems on the existence and convergence are proved.

\section{Problem Formulation}
\setcounter{equation}{0}
\label{formulation}
For the rate of convergence, we focus on the following Bolza problem of control systems in the feedback linearizable normal form. A more complicated problem with constraints is studied in Section \ref{sec-convergence} for the existence and convergence of approximate solutions. 
\vspace{0.1in}

\noindent{\bf Problem B:}\ \ 
Determine the state-control function pair $(x(t),u(t))$, $x\in \R^r$ and $u\in \R $, that minimizes the cost function
\bea
J(x(\cdot),u(\cdot))&  = & \int_{-1}^{1}F(x(t), u(t))\ dt + E(x(-1),x(1)) \label{eqad3a}
\eea
subject to the following differential equations and initial condition
\bea
&&\left\{ 
\begin{array}{lll}
\dot x_1  =  x_2\\
\;\;\; \vdots  \\
\dot x_{r-1}  = x_r\\
\dot x_r  =  f(x) + g(x)u 
\end{array}\label{eqad3}\right.\\
&& x(-1)= x_0\label{eqad3ini}
\eea
where $x\in \R^r$, $u\in \R $, and $F: \R^r \times \R \to \R$, $E: \R^r \times
\R^r \to \R$, $f: \R^r \to \R$, and $g: \R^r \to \R$ are all Lipschitz continuous functions with respect to their arguments. In addition, we assume $g(x)\neq 0$ for all $x$. 

Throughout the paper we make extensive use of Sobolev spaces, $W^{m,p}$, that
consists of functions, $\xi: [-1, 1] \to \mathbb{R}$ whose $j$-th order weak derivative, $\xi^{(j)}$, lies in $L^p$ for all $0 \le j \le m$
with the norm,
$$ \parallel \xi\parallel_{W^{m,p}} \quad = \sum_{j=0}^m \parallel \xi^{(j)}\parallel_{L^p}  $$ 
In this paper, we only consider the problems that have at least one optimal solution in which $x_r^\ast(t)$ has bounded $m$-th order weak derivative, i.e. $x_r^\ast(t)$ is in $W^{m,\infty}$. For some results, we assume $m\geq 3$. For others, $m$ is smaller. Unless the term `strong derivative' is emphasized, all derivatives in the paper are in the weak sense.  

\begin{figure}[ht!]
      \begin{center}
      \includegraphics[width=12cm] {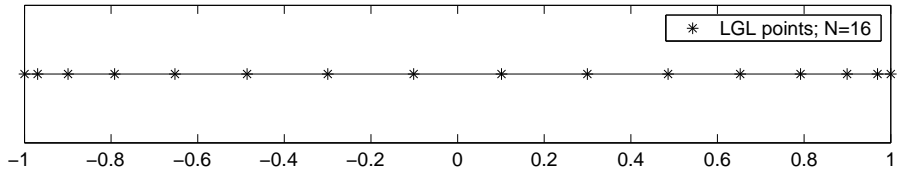}
      \caption{LGL nodes $N=16$ }      \label{figure-nodes}
      \end{center}
\end{figure}

The PS optimal control method addressed in this paper is an efficient direct method. In typical direct methods, the original optimal control problem, not the associated necessary conditions, is discretized to formulate a nonlinear programming problem. The accuracy of the discretization is largely determined by the accuracy of the underlying approximation method. Given any function $f(t): [a,b] \rightarrow \R$, a conventional method of approximation is to interpolate at uniformly spaced nodes: $t_0=a$, $t_1=(b-a)/N$, $\cdots$, $t_N=b$. However, it is known that uniform spacing is not efficient. More sophisticated node selection methods are able to achieve significantly improved accuracy with fewer nodes. It is important to emphasize that, for optimal control problems, the rate of convergence is not merely an issue of efficiency; more importantly it is about feasibility. An increased number of nodes in discretization results in a higher dimension in the nonlinear programming problem. A computational method becomes practically infeasible when the dimension and complexity of the nonlinear programming exceed the available computational power. In a PS approximation based on Legendre-Gauss-Lobatto (LGL) quadrature nodes, a function $f(t)$ is  approximated by $N$-th order Lagrange polynomials using the interpolation at these nodes. The LGL nodes, $t_0=-1 <t_1<\cdots < t_N=1$, are defined by
$$\begin{array}{llll}
t_0=-1,\;\; t_N=1, \mbox{ and }\\
\mbox{for }k=1,2,\ldots,N-1, t_k \mbox{ are the roots of } \dot L_N(t)
\end{array}
$$
where $\dot L_N(t)$ is the derivative of the $N$-th order Legendre polynomial $L_N(t)$. The discretization works in the interval of $[-1, 1]$. An example of LGL nodes with $N=16$ is shown in Figure \ref{figure-nodes}.
It was proved in approximation theory that the polynomial interpolation at the LGL nodes converges to $f(t)$ under $L^2$ norm at the rate of $1/N^{m}$, where $m$ is the smoothness of $f(t)$ (see for instance \cite{canuto} Section 5.4). If $f(t)$ is $C^\infty$, then the polynomial interpolation at the LGL nodes converges at a spectral rate, i.e. it is faster than any given polynomial rate. This is a very impressive convergence rate. 

PS methods have been widely applied in scientific computation for models governed by partial differential equations, such as complex fluid dynamics. However, PS optimal control has several fundamental differences from the computation of PDEs. Solving optimal control problems asks for the approximation of several objects collectively, including the differential equation that defines the control system, the integration in the cost function, and the state and control trajectories. In addition to the various types of approximations, a nonlinear programming must be applied to the overall discretized optimization problem to find an approximate optimal control. All these factors may deteriorate the final approximate solution. The existing theory of PS approximation of differential equations is not applicable. New theory needs to be developed for the existence, convergence, and the rate of convergence for optimal control problems.  

In the following, we introduce the notations used in this paper. Then, the discretized nonlinear programming problem is formulated. In a PS optimal control method, the state and control functions, $x(t)$ and $u(t)$, are approximated by $N$-th order Lagrange polynomials based on the interpolation at the LGL quadrature nodes. In the discretization, the state variables are approximated by the vectors $\bar x\NK \in \R^r$, i.e.
$$\bar x\NK=\MT \bar x_1\NK\\ \bar x_2\NK\\ \vdots \\ \bar x_r\NK\EM$$
is an approximation of $x(t_k)$. Similarly, $\bar u\NK$ is the approximation of $u(t_k)$. Thus, a discrete approximation of the function $x_i(t)$ is the vector 
$$\bar x_i^N=\MT \bar x_i^{N1}& \bar x_i^{N2} & \cdots & \bar x_i^{NN}\EM$$ 
A continuous approximation is defined by its polynomial interpolation, denoted by $x_i^N(t)$, i.e.
\bea
x_i(t)& \approx & x_i^N(t)=\sum_{k=0}^{N}\bar x_i^{Nk}\phi_k(t), \label{eqad6a}
\eea
where $\phi_k(t)$ is the Lagrange interpolating polynomial \cite{canuto}.
Instead of polynomial interpolation, the control input is approximated by the following non-polynomial interpolation
\bea
\label{approxu}
u^N(t)=\Fr{\dot x_r^N(t)-f(x^N(t))}{g(x^N(t))}
\eea
In the notations, the discrete variables are denoted by letters with an upper bar, such as $\bar x^{Nk}_i$ and $\bar u^{Nk}$. If $k$ in the superscript and/or $i$ in the subscript are missing, it represents the corresponding vector or matrix in which the indices run from minimum to maximum. For example,
\bean
\bar x^{N}_i&=&\MT \bar x_i^{N0} & \bx_i^{N1} & \cdots & \bx_i^{NN}\EM\\
\bar x^{Nk}&=& \MT \bx_1^{Nk}\\ \bx_2^{Nk}\\ \vdots \\ \bx_r^{Nk}\EM\\
\bar x^N &=& \MT \bx_1^{N0} & \bx_1^{N1} & \cdots & \bx_1^{NN}\\
\bx_2^{N0} & \bx_2^{N1} & \cdots & \bx_2^{NN}\\
\vdots & \vdots &\vdots & \vdots \\
\bx_r^{N0} & \bx_r^{N1} & \cdots & \bx_r^{NN}\EM
\eean
Similarly,
$$\bar u^{N}=\MT \bu^{N0}& \bu^{N1} & \cdots &\bu^{NN}\EM$$
Given a discrete approximation of a continuous function, the interpolation is denoted by the same notation without the upper bar. For example, $x_i^N(t)$ in (\ref{eqad6a}), $u^N(t)$ in (\ref{approxu}). The superscript $N$ represents the number of LGL nodes used in the approximation. Throughout this paper, the interpolation of $(\bx^N,\bu^N)$ is defined by (\ref{eqad6a})-(\ref{approxu}), in which $u^N(t)$ is not necessarily a polynomial. It is proved in Lemma \ref{lemma1} that (\ref{approxu}) is indeed an interpolation. 

Existing results in the analysis of spectral methods show that PS method is an approach that is easy and accurate in the approximation of smooth functions, integrations, and differentiations, all critical to optimal control problems. For differentiation, the derivative of $x^N_i(t)$ at the LGL node $t_k$ is easily computed by the following matrix multiplication \cite{canuto}
\bea
\MT \dot x_i^N(t_0) & \dot x_i^N(t_1) & \cdots & \dot x_i^N(t_N)\EM^T=D(\bx^N_i)^T \label{eqDxN}
\eea
where the $(N+1)\times(N+1)$ differentiation matrix $D$ is defined by
\bean
D_{ik} & = & \left\{\begin{array}{ll} \frac{L_N(t_i)}{L_N(t_k)}\frac{1}{t_i-t_k},&  \mbox{if}\ \ i\neq k; \\
                                       \\
                                       -\frac{N(N+1)}{4}, & \mbox{if}\ \ i=k=0;\\
                                       \\
                                       \frac{N(N+1)}{4},& \mbox{if}\ \ i=k=N;\\
                                       \\
                                       0, & \mbox{otherwise} \end{array} \right.
\eean
The cost functional $J[x(\cdot), u(\cdot)]$ is approximated by the
Gauss-Lobatto integration rule,
\bea J[x(\cdot),u(\cdot)] \ \approx \ \bar J^N(\bx^N,\bu^N) & = &
\sum_{k=0}^{N} F(\bx^{Nk},\bu^{Nk})w_k + E(\bx^{N0},\bx^{NN}) \nn \eea
where $w_k$ are the LGL weights defined by \bea w_k & = &
\frac{2}{N(N+1)}\frac{1}{[L_N(t_k)]^2}, \nn \eea 
The approximation is so accurate that it has zero error if the integrand function is a polynomial of degree less than or equal to $2N-1$, a degree that is almost a double of the number of nodes \cite{canuto}. Now, we are ready to define Problem \BN, a PS discretization of Problem B. 

For any integer $m_1 >0$, let $\{a_0^{N}(m_1),a_1^{N}(m_1), \cdots, a_{N-r-m_1+1}^{N}(m_1)\}$ denote the coefficients in the Legendre polynomial expansion for the interpolation polynomial of the vector $\bar x_r^N(D^T)^{m_1}$. Note that the interpolation of $\bar x_r^N(D^T)^{m_1}$ equals the polynomial of $\frac{d^{m_1}x_r^N(t)}{dt^{m_1}}$. Thus, there are only $N-r-m_1+2$ nonzero spectral coefficients because it is proved in Section \ref{rate} that the order of $\frac{d^{m_1} x^N_r(t)}{dt^{m_1}}$ is at most degree of $N-r-m_1+1$. These coefficients depend linearly on $\bx_r^N$ \cite{boyd}, 
\EQ
{\tiny
\MT
  a^{N}_0(m_1) \\
  \vdots \\
  a^{N}_{N-r-m_1+1}(m_1)
\EM 
= }\\
{\tiny
\MT
  \frac12 &  &   \\
   & \ddots &   \\   
   & & N-r-m_1+1+\frac12 
\EM
\MT  L_0(t_0) & \cdots & L_0(t_N) \\
   & \vdots &  \\
  L_{N-r-m_1+1}(t_0) & \cdots & L_{N-r-m_1+1}(t_N)
\EM
\MT  w_0 &  &   \\
   & \ddots &   \\   
   & & w_N \\
\EM D^{m_1}
\MT   \bar x^{N0}_{r} \\
  \vdots \\
  \bar x^{NN}_{r} \\
\EM
}
\label{spectralcoef}
\EE
The PS discretization of Problem \BN is defined as follows.
\vspace{0.1in}

\noindent{\bf Problem ${\bf B}^{\bf N}$:}\ \ 
Find $\bar x\NK\in \R^{r}$ and $\bar u\NK\in \R$, $k \ = \ 0,1,\ldots,N$, that minimize
\bea \bar J^N(\bar x^N ,\bar u^N) & = & \sum_{k=0}^{N}F(\bar x\NK,\bar u\NK)w_k + 
E(\bar x^{N0},\bar x^{NN}) \label{JN1}
\eea
subject to
\bea
&&\left\{ \begin{array}{rcl}
D ( \bar x_1^N)^T &=& (\bar x_2^N)^T \\
D (\bar x_2^N)^T&=&(\bar x_3^N)^T \\
& \vdots &\\
D (\bar x_{r-1}^N)^T&=&(\bar x_r^N)^T \\
D (\bar x_r^N)^T& = & \MT f(\bar x^{N0})+g(\bar x^{N0})\bar u^{N0}\\ \vdots \\
  f(\bar x^{NN})+g(\bar x^{NN})\bar u^{NN} \EM \\
  \end{array}\right.\label{eq3_2}\\
  \nn\\
&&\bar x^{N0}=x_0\label{eq3_3}\\
\nn\\
 &&\underline{\bb}  \leq  \MT \bar x\NK \\ \bar u\NK\EM
\ \leq \ \bar \bb, \;\; \;\; \mbox{ for all } 0\leq k\leq N\label{eq3_4b}\\
&&\underline{\bb}_j \leq  \MT 1&0& \cdots &0\EM D^j (\bx_r^N)^T \ \leq \bar \bb_j, \mbox{ if }1\leq j\leq m_1-1 \mbox{ and } m_1\geq 2 \label{eq3_4a}\\
&& \dsum_{n=0}^{N-r-m_1+1}|a^{N}_n(m_1) |\leq \bd \label{eq3_4}
\eea
\vspace{0.1in}

Comparing to Problem B, (\ref{eq3_2}) is the discretization of the control system defined by the differential equation. The regularization condition (\ref{eq3_4}) assures that the derivative of the interpolation up to the order of $m_1$ is bounded. It is proved in the following sections that the integer $m_1$ is closely related to the convergence rate. The inequalities (\ref{eq3_4b}), (\ref{eq3_4a}) and (\ref{eq3_4}) are regularization conditions that do not exist in Problem B. It is proved in the next few sections that these additional constraints do not affect the feasibility of Problem \BN . Therefore, it does not put an extra limit to the family of problems to be solved. 

In searching for a discrete optimal solution, it is standard for software packages of nonlinear programming to  require a search region. Typically, the search region is defined by a constraint (\ref{eq3_4b}). However, this box-shaped region may contain solutions that are not good approximations of the continuous-time solution. To guarantee the rate of convergence, the search region is refined to a smaller one by imposing constraints (\ref{eq3_4a}) and (\ref{eq3_4}). It is proved in this paper that there always exist feasible solutions that satisfy all the constraints and the optimal cost converges, provided the upper and lower bounds are large enough. In (\ref{eq3_4a}), $\underline{\bb}_j$ and $\bar \bb_j$ represent the bounds of initial derivatives. In (\ref{eq3_4}), $\bd $ is the bound determined by $x_r^{(m_1+1)}$ satisfying the inequality (\ref{eqdbound}). Without known the optimal solution, these bounds of search region have to be estimated before computation or they are determined by numerical experimentations. The constraints (\ref{eq3_4a}) and (\ref{eq3_4}) are necessary to avoid the restrictive consistent approximation assumption made in \cite{gong}. At a more fundamental level, the order of derivatives, $m_1$ in (\ref{eq3_4}), determines the convergence rate of the approximate optimal control. Another interesting fact that amply justify these additional constraints is that Problem \BN\ may not even have an optimal solution if we do not enforce (\ref{eq3_4b}). This is shown by the following counter example.

\begin{example} \rm 
Consider the following problem of optimal control.
\bea
&&\min_{(x(\cdot),u(\cdot))} \int_{-1}^1 \Fr{(x(t)-u(t))^2}{u(t)^4}dt\nn\\
&&\dot x=u\label{eqex}\\
&&x(-1)=e^{-1}\nn
\eea
It is easy to check that the optimal solution is 
\bean
u=e^t, & x(t)=e^t
\eean
and the optimal cost value is zero. Although the solution to the problem (\ref{eqex}) is simple and analytic, the PS discretization of (\ref{eqex}) does not have an optimal solution if the constraint (\ref{eq3_4b}) is not enforced. To prove this claim, consider the PS discretization,
\bea
&&\min_{(\bx^N, \bu^N)} \bar J^N(\bx^N,\bu^N)=\dsum_{k=0}^{N} \Fr{(\bx^{Nk}-D_k (\bx^N)^T)^2}{\left( D_k(\bx^N)^T\right)^4}w_k\nn\\
&&D(\bx^N)^T=(\bu^N)^T\label{eqex1}\\
&&\bx^{N0} =e^{-1}\nn
\eea
where $D_k$ is the $k$th row of the differentiation matrix $D$. Let $x^N(t)$ be the interpolation polynomial of $\bx^N$, then it is obvious that 
$$x^N(t)-\dot x^N(t)\not\equiv 0$$
Thus, there exists $k$ so that 
$$\bx^{Nk}-D_k (\bx^N)^T \neq 0$$
So, 
\EQ
\label{eqex2}
\bar J^N(\bx^N,\bu^N)>0
\EE
for all feasible pairs $(\bx^N, \bu^N)$. For any $\alpha>0$, define
$$\bx\NK = e^{-1} +\alpha (t_k+1)$$
The interpolation of $\bx^N$ is the linear polynomial 
$$x^N(t)=e^{-1}+\alpha (t+1)$$
Then, 
$$D_k (\bx^N)^T =\dot x^N(t_k)=\alpha$$
The cost function is
\bean
\bar J^N(\bx^N,\bu^N)&=&\dsum_{k=0}^N \Fr{(e^{-1}+\alpha (t_k+1)-\alpha)^2}{\alpha^4}w_k\\
&=& \dsum_{k=0}^N\Fr{(e^{-1}+\alpha t_k)^2}{\alpha^4}w_k\\
&\leq & \dsum_{k=0}^N\Fr{(e^{-1}+\alpha)^2}{\alpha^4}w_k\\
&=& 2\Fr{(e^{-1}+\alpha)^2}{\alpha^4}
\eean
Therefore, $\bar J^N(\bx^N,\bu^N)$ can be arbitrarily small as $\alpha$ approaches $\infty$. However, $\bar J^N(\bx^N,\bu^N)$ is always positive as shown by (\ref{eqex2}). We conclude that the discretization (\ref{eqex1}) has no minimum value for $\bar J^N(\bx^N,\bu^N)$. 
\end{example} 

\section{Convergence Rate}
\setcounter{equation}{0}
\label{rate}

Given a solution to Problem \BN, we use (\ref{approxu}) to approximate the optimal control. In this section we prove that, under this approximate optimal control, the value of the cost function converges to the optimal cost of Problem B as the number of nodes is increased. More importantly, we can prove a high-order rate of convergence. In the literature, it has been proved that PS methods have a spectral rate when approximating $C^\infty$ functions, i.e. the rate is faster than any polynomial rate. However, there are no results in the literature thus far on the convergence rate of PS optimal control. Meanwhile, in many problems solved by PS optimal control we clearly observed a rate of high-order in the convergence. In this section, we prove a convergence rate that depends on the smoothness of the optimal control. More specifically, the rate is about $\frac{1}{N^{2m/3-1}}$, where $m$ is defined by the smoothness of the optimal trajectory. If the cost function can be accurately computed, then the convergence rate is improved to $\frac{1}{N^{2m-1}}$. In the special case of $C^\infty$, it is proved that PS method is able to converge faster than any given polynomial rate. Before we introduce the main theorems of this section, the following example in \cite{gong} is briefly presented to show the rapid convergence of the PS optimal control method. 

\begin{example} \label{example-comp}\ \ \rm
Consider the following nonlinear optimal control problem:
\begin{eqnarray}\nn
&\left\{\begin{array}{lrl} {\rm Minimize } & J[x(\cdot), u(\cdot)]
=&4x_1(2)+x_2(2)+4\displaystyle \int_0^2 u^2(t) \ dt  \\
{\rm Subject\ to} & \dot x_1(t) = & x_2^3(t)\\
& \dot x_2(t) = & u(t)\\
&(x_1(0),x_2(0))= & (0,1)
\end{array}\right.&
\end{eqnarray}
It can be shown that the exact optimal control is defined by
$u^\ast(t) = - \frac{8}{(2+t)^3}$. For this problem, the PS method achieves the accuracy in the magnitude of $10^{-8}$ with only 18 nodes. A detail comparison of the PS method with some other discretization methods are addressed in \cite{gong}. From  Figure \ref{figure-exp1} in logarithmically scaled coordinates, it is obvious that the computation using the PS method converges exponentially.
\begin{figure}[h]
      \begin{center}
      \includegraphics[width=12cm] {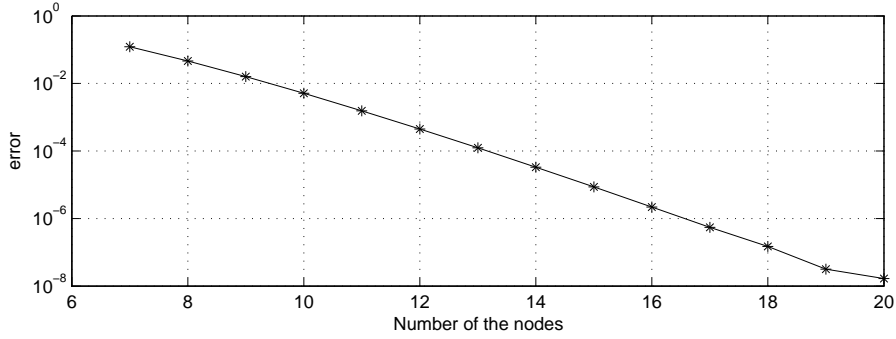}
      \caption{Error vs number of the nodes for the pseudospectral method }      \label{figure-exp1}
      \end{center}
\end{figure}
\end{example}

Problem \BN\ has several bounds in its definition, $\underline{\bb}$, $\bar \bb$, $\underline{\bb}_j$, $\bar \bb_j$, and $\bd$. These bounds can be selected from a range determined by Problem B. The constraints $\underline{\bb}$ and $\bar \bb$ are lower and upper bounds so that the optimal trajectory of Problem B is contained in the interior of the region. Suppose Problem B has an optimal solution $(x^\ast (t), u^\ast (t))$ in which $(x_r^{\ast}(t))^{(m)}$ has bounded variation for some $m\geq 3$, where $x_r^{\ast}(t)$ is the $r$th component of the optimal trajectory. Suppose $m_1$ in Problem \BN\ satisfies $2\leq m_1\leq m-1$. Then, we can select the bounds $\underline{\bb}_j$ and $\bar \bb_j$ so that $(x_r^\ast (t))^{(j)}$ is contained in the interior of the region. For $\bd$, we assume
\bea
\bd &>&\Fr{6}{\sqrt{\pi}}(U(\left. x_r^{\ast}\right.^{(m_1+1)})+V(\left.x_r^{\ast}\right.^{(m_1+1)}))\zeta(3/2)\label{eqdbound}
\eea
where $U(\left. x_r^{\ast}\right.^{(m_1+1)})$ is the upper bound and $V(\left.x_r^{\ast}\right.^{(m_1+1)})$ is the total variation of $\left. x_r^{\ast}\right.^{(m_1+1)}(t)$; and $\zeta(s)$ is the $\zeta$ function defined by
\bea
\label{zetafun}
&&\zeta (s)=\sum_{k=1}^\infty \Fr{1}{k^s}
\eea
If all the bounds are selected as above, then it is proved in Section \ref{sec-convergence} that Problem \BN\ is always feasible provided $m\geq 2$. Note that in practical computation, $\underline{\bb}$, $\bar \bb$, $\underline{\bb}_j$, $\bar \bb_j$, and $\bd$ are unknown. They must be estimated based upon experience or other information about the system. 

\begin{theorem}
\label{th_rate}
Suppose Problem B has an optimal solution $(x^\ast (t), u^\ast (t))$ in which the $m$-th order derivative $(x_r^{\ast}(t))^{(m)}$ has a bounded variation for some $m\geq 3$. In Problem \BN , select $m_1$ and $\alpha$ so that $1\leq m_1\leq m-1$ and $0< \alpha <m_1-1$. Suppose $f(\cdot)$, $g(\cdot)$, $F(\cdot)$, and $E(x_0,\cdot)$ are $C^{m}$ and globally Lipschitz. Suppose all other bounds in Problem \BN\ are large enough. Given any sequence
\EQ
\label{eq3_13}
\{ (\bar x^{\ast N},\bar u^{\ast N})\}_{N\geq N_1}
\EE
of optimal solutions of Problem \BN . Then the approximate cost converge to the optimal value at the following rate
\bea
\left| J(x^\ast(\cdot),u^\ast (\cdot))-J( x^{\ast N}(\cdot), u^{\ast N}(\cdot))\right|
&\leq& \Fr{M_1}{(N-r-m_1-1)^{2m-2m_1-1}}+\Fr{M_2}{N^{\alpha}}\label{eq3_12b}\\
\nn\\
\left| J(x^\ast(\cdot),u^\ast (\cdot))-\bar J^{N}(\bar x^{\ast N},\bar u^{\ast N})\right|
&\leq& \Fr{M_1}{(N-r-m_1-1)^{2m-2m_1-1}}+\Fr{M_2}{N^{\alpha}}\label{eq3_12}
\eea
where $M_1$ and $M_2$ are some constants independent of $N$. In (\ref{eq3_12b}), $(x^{\ast N}(t), u^{\ast N}(t))$ is the interpolation of (\ref{eq3_13}) defined by (\ref{eqad6a})-(\ref{approxu}). In fact, $x^{\ast N}(t)$ is the trajectory of (\ref{eqad3}) under the control input $u^{\ast N}(t)$.
\end{theorem}

Theorem \ref{th_rate} implies that the costs of any sequence of discrete optimal solutions must converge to the optimal cost of Problem B, no matter the sequence of the discrete state and control trajectories converge or not. In other words, it is possible that the sequence of discrete optimal controls does not converge to a unique continuous-time control; meanwhile the costs using these approximate optimal controls converge to the true optimal cost of Problem B. Therefore, this theorem does not require the local uniqueness of solutions for Problem B. This is different from many existing convergence theorems of computational optimal control, in which a unique optimal solution and coercivity are assumed. This is made possible because the proofs in this paper do not rely on the necessary conditions of optimal control. The key idea in the proof is to shape the search region in Problem \BN by regulating the discrete solutions using (\ref{eq3_4b})-(\ref{eq3_4a})-(\ref{eq3_4}). We would like to emphasize that the regulation constraints are added to the discretized problem, not the original Problem B. So, the constraints do not restrict the problem to be solved, and they do not need to be verified before computation. In addition, increasing the number of constraints results in smaller search region for an optimal solution.

\begin{remark}
If $f(\cdot)$, $g(\cdot)$, $F(\cdot)$ and $x_r^\ast(t)$ are $C^\infty$, then we can select $m$ and $m_1$ arbitrarily large. In this case, we can make the optimal cost of Problem \BN converge faster than any given polynomial rate. 
\end{remark}

\begin{remark}
From (\ref{eq3_12b}) and (\ref{eq3_12}), the convergence rate is determined by $m$, the smoothness of the optimal trajectory, and $m_1$, the order in the regulation of discrete solutions. While $m$ is a property of Problem B that cannot be changed, $m_1$ in Problem \BN\ can be selected within a range. However, the errors in (\ref{eq3_12b}) and (\ref{eq3_12}) have two parts, one is an increasing function of $m_1$ and the other is a decreasing function of $m_1$. In Corollary \ref{corollary1}, we show an optimal selection of $m_1$ to maximize the combined convergence rate. 
\end{remark}

The proof is convoluted involving results from several different areas, including nonlinear functional analysis, orthogonal polynomials, and approximation theory. First, we introduce the concept of Fr\'echet derivative. Let us consider the continuous cost function, $J(x(\cdot), u(\cdot))$, subject to (\ref{eqad3})-(\ref{eqad3ini}) as a nonlinear functional of $u(\cdot)$, denoted by ${\cal J}(u)$. For any $u$ in the Banach space $W^{m-1,\infty}$, suppose there exists a linear bounded operator ${\cal L}$: $W^{m-1,\infty}\rightarrow \R$ such that 
$$|{\cal J}(u+\Delta u)-{\cal J}(u) - L \Delta u|=o(||\Delta u||_{W^{m-1,\infty}})$$
for all $u+\Delta u$ in an open neighborhood of $u$ in $W^{m-1,\infty}$. Then, ${\cal L}$ is called the Fr\'echet derivative of ${\cal J}(u)$ at $u$, denoted by ${\cal J}^\prime (u)=\cal L$. If ${\cal J}^\prime (u)$ exists at all points in an open subset of $W^{m-1,\infty}$, then 
${\cal J}^\prime (u)$ is a functional from this open set to the Banach space $L(W^{m-1,\infty}, \R)$ of all bounded linear operators. If this new functional has a Fr\'echet derivative, then it is called the second order Fr\'echet derivative, denoted by ${\cal J}^{\prime\prime} (u)$. The following lemma is standard in nonlinear functional analysis \cite{wouk}. 

\begin{lemma}
\label{lemmafrechet}
Suppose ${\cal J}$ takes a local minimum value at $u^\ast$. Suppose ${\cal J}$ has second order Fr\'echet derivative at $u^\ast$. Then, 
$${\cal J}(u^\ast+\Delta u)=({\cal J}^{\prime\prime} (u^\ast)\Delta u)\Delta u +o(||\Delta u||^2)$$
\end{lemma}

The rate of convergence for the spectral coefficients can be estimated by the following Jackson's Theorem.
\begin{lemma}
\label{Jackson}
(Jackson's Theorem \cite{sansone}) Let $h(t)$ be of bounded variation in $[-1,1]$. Define
$$H(t)=H(-1)+\displaystyle{\int_{-1}^t}h(s)ds$$
then $\{ a_n \}_{n=0}^\infty$, the sequence of spectral coefficients of $H(t)$, satisfies the following inequality
$$a_n<\Fr{6}{\sqrt{\pi}}(U(h(t))+V(h(t)))\Fr{1}{n^{3/2}}$$
for $n\geq 1$.
\end{lemma}

Given a continuous function $h(t)$ defined on $[-1, 1]$. Let $\hat p^N(t)$ be the best polynomial of degree $N$, i.e. the $N$th order polynomial with the smallest distance to $h(t)$ under $||\cdot ||_\infty$ norm. Let $I_Nh(t)$ be the polynomial interpolation using the value of $h(t)$ at the LGL nodes. Then, we have the following inequality from the theory of approximation and orthogonal polynomials \cite{canuto}, \cite{hesthaven}.
\begin{lemma}
\label{lebesgue}
$$|| h(t)-I_Nh||_\infty \leq (1+\Lambda_N)||h(t) - \hat p^N(t)||_\infty$$
where $\Lambda_N$ is called Lebesgue constant. It satisfies
$$
\Lambda_N \leq \Fr{2}{\pi}log(N+1)+0.685\cdots
$$
\end{lemma}

The best polynomial approximation represents the closest polynomial to a function under $||\cdot ||_\infty$. The error can be estimated by the following Lemma \cite{canuto}.

\begin{lemma} 
\label{best}
(1) Suppose $h(t)\in W^{m,\infty}$. Let $\hat p^N(t)$ be the best polynomial approximation. Then
$$||\hat p^N(t)-h(t)||_\infty\leq \Fr{C}{N^{m}} ||h(t)||_{W^{m, \infty}}$$
for some constant $C$ independent of $h(t)$, $m$ and $N$. \\
\\
(2) If $h(t)\in W^{m,2}$, then 
$$|| h(t)-P_N h (t)||_\infty \leq \Fr{C || h(t)||_{W^{m,2}}}{N^{m-3/4}}$$
where $P_Nh$ is the N-th order truncation of the Legendre series of $h(t)$. \\
(3) If $h(t)$ has the $m$-th order strong derivative with a bounded variation, then 
$$|| h(t)-P_N h (t)||_\infty \leq \Fr{C V(h^{(m)}(t))}{N^{m-1/2}}$$
\end{lemma}

The following lemmas are proved specifically for PS optimal control methods. Similar results can be found in \cite{kang} except that some assumptions on $m_1$ are relaxed.   

\begin{lemma}
\label{lemma1} (\cite{kang}) (i) For any trajectory, $(\bar x^N,\bar u^N)$, of the dynamics (\ref{eq3_2}), the pair $(x^N(t), u^N(t))$ defined by (\ref{eqad6a})-(\ref{approxu}) satisfies the differential equations defined in (\ref{eqad3}). Furthermore, 
\bea
&&\bx^{Nk}=x^N(t_k), \; \bu^{Nk}=u^N(t_k), \; \mbox{for } k=0,1,\cdots,N\label{eqmap}
\eea  

(ii) For any pair $(x^N(t),u^N(t))$ in which $x^N(t)$ consists of polynomials of degree less than or equal to $N$ and $u^N(t)$ is a function, if $(x^N(t),u^N(t))$ satisfies the differential equations in (\ref{eqad3}), then $(\bx^N,\bu^N)$ defined by (\ref{eqmap}) satisfies (\ref{eq3_2}).

(iii) If $(\bar x^N,\bar u^N)$ satisfies (\ref{eq3_2}), then the degree of $x_i^N(t)$ is less than or equal to $N-i+1$.
\end{lemma}
\begin{proof} (i) Suppose $(\bar x^N,\bar u^N)$ satisfies the equations in (\ref{eq3_2}). Because $x^N(t)$ is the polynomial interpolation of $\bx^N$, and because of equations (\ref{eqDxN}), we have
\EQ
\MT \dot x_i^N(t_0) & \dot x_i^N(t_1)&\cdots & \dot x_i^N(t_N)\EM\\
=\bar x_i^N D^T \\
=\bar x_{i+1}^N\\
=\MT x_{i+1}^N(t_0) & x_{i+1}^N(t_1)&\cdots & x_{i+1}^N(t_N)\EM\nn
\EE
Therefore, the polynomials $\dot x^N_i(t)$ and $x^N_{i+1}(t)$ must equal each other because they coincide at $N+1$ points and because the degrees of $x_i^N(t)$ and $x_{i+1}^N(t)$ are less than or equal to $N$. In addition, (\ref{approxu}), the definition of $u^N(t)$, implies the last equation in (\ref{eqad3}). So, the pair $(x^N(t),u^N(t))$ satisfies all equations in (\ref{eqad3}). Now, we prove (\ref{eqmap}). Because $x^N(t)$ is an interpolation of $\bx^N$, we know  $\bx^{Nk}=x^N(t_k)$ for $0\leq k\leq N$. From (\ref{approxu}), 
\bea
u^N(t_k)&=&\Fr{\dot x_r^N(t_k)-f(x^N(t_k))}{g(x^N(t_k))}\nn\\
&=& \Fr{\dot x^N_r(t_k)-f(\bar x\NK)}{g(\bar x\NK)}\label{eqadstar}
\eea
Because of (\ref{eqDxN}), we have
$$\MT \dot x_r^N(t_0) & \dot x_r^N(t_1)&\cdots & \dot x_r^N(t_N)\EM^T= D (\bar x_r^N)^T$$
Therefore, (\ref{eqadstar}) is equivalent to
\bean
\MT u^N(t_0) & u^N(t_1)&\cdots & u^N(t_N)\EM^T &=& \mbox{diag}\left(\Fr{1}{g(\bar x^{N0})}, \cdots, \Fr{1}{g(\bar x^{NN})}\right)\left( D(\bar x_r^N)^T-\MT f(\bar x^{N0})\\ \vdots \\
  f(\bar x^{NN})\EM\right)
\eean
Comparing to the last equation in (\ref{eq3_2}), it is obvious that $u^N(t_k)=\bar u\NK$. So, (\ref{eqmap}) holds true. Part (i) is proved. 

(ii) Assume $(x^N(t),u^N(t))$ satisfies the differential equations in (\ref{eqad3}). Because $x^N(t)$ are polynomials, (\ref{eqDxN}) implies  
\bea
&&\bar x_i^N D^T \nn\\
&&=\MT \dot x_i^N(t_0) & \dot x_i^N(t_1)&\cdots & \dot x_i^N(t_N)\EM\label{eqad1}\\
&&=\MT x_{i+1}^N(t_0) & x_{i+1}^N(t_1)&\cdots & x_{i+1}^N(t_N)\EM\nn\\
&&=\bar x_{i+1}^N\nn
\eea
Furthermore, 
\bea
\bx^N_r D^T&=&\MT \dot x_r^N(t_0)&\dot x_r^N(t_1)&\cdots &\dot x_r^N(t_N)\EM\label{eqad2}\\
&=& \MT f(x^N(t_0))+g(x^N(t_0))u^N(t_0)&\cdots & f(x^N(t_N))+g(x^N(t_N))u^N(t_N)\EM\nn
\eea
Equations (\ref{eqad1}) and (\ref{eqad2}) imply that $(\bx^N,\bu^N)$ satisfies (\ref{eq3_2}).
Part (ii) is proved.

(iii) We know that the degree of $x_1^N(t)$, the interpolation polynomial, is less than or equal to $N$. From (i), we know $x_2^N(t)=\dot x_1^N(t)$. Therefore, the degree of $x_2^N(t)$ must be less than or equal to $N-1$. In general, the degree of $x_i^N(t)$ is less than or equal to $N-i+1$. 
\end{proof}

\begin{lemma}
\label{lemma1.5}
Suppose $\{ (\bx^N,\bu^N)\}_{N=N_1}^\infty$ is a sequence satisfying (\ref{eq3_2}), (\ref{eq3_4b}), (\ref{eq3_4a}) and (\ref{eq3_4}), where $m_1\geq 1$. Then,
\bean
&&\left\{ \left. || (x^N(t))^{(l)}||_\infty \right| N\geq N_1,\, l=0,1,\cdots,m_1\right\}
\eean 
is bounded. If $f(x)$ and $g(x)$ are $C^{m_1-1}$, then
\bean
&&\left\{ \left. || (u^N(t))^{(l)}||_\infty \right| N\geq N_1,\, l=0,1,\cdots,m_1-1\right\}
\eean 
is bounded.
\end{lemma}

\begin{proof}
Consider $(x_r^N(t))^{(m_1)}$. From Lemma \ref{lemma1}, it is a polynomial of degree less than or equal to $N-r-m_1+1$. Therefore, 
\bean
&&(x_r^N(t))^{(m_1)} = \dsum_{n=0}^{N-r-m_1+1} a_n^{N}(m_1)L_n(t)
\eean
where $L_n(t)$ is the Legendre polynomial of degree $n$. It is known that $|L_n(t)|\leq 1$. Therefore, (\ref{eq3_4}) implies that $||(x_r^N(t))^{(m_1)}||_\infty$ is bounded by $\bd$ for all $N\geq N_1$. Now, let us consider $(x_r^N(t))^{(m_1-1)}$. From (\ref{eqDxN}) we have,
\bean
(x_r^N(t))^{(m_1-1)}&=&(x_r^N(t))^{(m_1-1)}|_{t=-1}+\int_0^t (x_r^N(s))^{(m_1)}ds\\
&=&\MT 1& 0&\cdots &0\EM D^{m_1-1}(\bx_r^N)^T + \int_0^t (x_r^N(s))^{(m_1)}ds
\eean
So, $||(x_r^N(t))^{(m_1-1)}||_\infty$, $N\geq N_1$, is bounded because of (\ref{eq3_4a}). Similarly, we can prove all derivatives of $x_r^N(t)$ of order less than $m_1$ are bounded. The same approach can also be applied to prove the bound
$$u^N(t)=\Fr{\dot x_r^N(t)-f(x^N(t))}{g(x^N(t))}$$
Because $f(x)$ and $g(x)$ have continuous derivatives of order less than or equal to $m_1-1$, the boundedness of 
$$\left\{ \left. || (u^N(t))^{(l)}||_\infty \right| N\geq N_1,\, j=0,1,\cdots,m_1-1\right\}$$ follows the boundedness of $(x_r^N(t))^{(l)}$ proved above.
\end{proof}

Given any function $h(t)$ defined on $[-1, 1]$. In the following, $U(h)$ represents an upper bound of $h(t)$ and  $V(h)$ represents the total variation.

\begin{lemma}
\label{lemmafeasibility}
Let $(x(t),u(t))$ be a solution of the differential equation (\ref{eqad3}). Suppose $x_r^{(m)}(t)$ has bounded variation for some $m\geq 2$. Let $m_1$ be an integer satisfying $1\leq m_1\leq m-1$. Then, there exist constants $M>0$ and $N_1>0$ so that for each integer $N\geq N_1$ the differential equation (\ref{eqad3}) has a solution $(x^N(t), u^N(t))$ in which $x^N(t)$ consists of polynomials of degree less than or equal to $N$. Furthermore, the pair $(x^N(t), u^N(t))$ satisfies 
\bea
||x_i^N(t)-x_i(t)||_\infty &\leq& \Fr{M|| x_r||_{W^{m,2}}}{(N-r-m_1+1)^{(m-m_1)-3/4}},  \;\;\;i=1,2,\cdots,r\label{eq1a}\\
||(x_r^N(t))^{(l)}-(x_r(t))^{(l)}||_\infty &\leq & \Fr{M|| x_r||_{W^{m,2}}}{(N-r-m_1+1)^{(m-m_1)-3/4}}, \;\;\;\; l=1,2,\cdots, m_1 \label{eq1d}\\
||u^N(t)-u(t)||_\infty &\leq&     \Fr{M|| x_r||_{W^{m,2}}}{(N-r-m_1+1)^{(m-m_1)-3/4}}\label{eq1db}
\eea
Furthermore, the spectral coefficients of $(x^N_r)^{(m_1)}(t)$ satisfy 
\bea
\label{eq1c}
|a^{N}_n(m_1)|\leq \Fr{6(U(x^{(m_1+1)}_r)+V(x^{(m_1+1)}_r))}{\sqrt{\pi} n^{3/2}}, \;\; n=1,2,\cdots,N-r-1
\eea
If $f(x)$ and $g(x)$ have Lipschitz continuous $L$th order partial derivatives for some $L\leq m_1-1$, then
\bea
||(u^N(t))^{(l)}-(u(t))^{(l)}||_\infty &\leq&     \Fr{M|| x_r||_{W^{m,2}}}{(N-r-m_1+1)^{(m-m_1)-3/4}}, \;\;\; l=1,\cdots,L\label{eq1b}
\eea
Furthermore, 
\bea
\label{eq1b2}
\begin{array}{lll}
x^N(-1)=x(-1)\\
u^N(-1)=u(-1),  &\mbox{ If } m_1\geq 2
\end{array}
\eea
\end{lemma}

\begin{remark}
\label{rem1}
In this lemma, if $x_r(t)$ has the $m$-th order strong derivative and if $x_r^{(m)}(t)$ has bounded variation for some $m\geq 2$, then the inequalities (\ref{eq1a}), (\ref{eq1d}), and (\ref{eq1db}) are slightly tighter.
\bea
||x_i^N(t)-x_i(t)||_\infty &\leq& \Fr{M|| x_r||_{W^{m,2}}}{(N-r-m_1+1)^{(m-m_1)-1/2}},  \;\;\;i=1,2,\cdots,r\label{eq1at}\\
||(x_r^N(t))^{(l)}-(x_r(t))^{(l)}||_\infty &\leq & \Fr{M|| x_r||_{W^{m,2}}}{(N-r-m_1+1)^{(m-m_1)-1/2}}, \;\;\;\; l=1,2,\cdots, m_1 \label{eq1dt}\\
||u^N(t)-u(t)||_\infty &\leq&     \Fr{M|| x_r||_{W^{m,2}}}{(N-r-m_1+1)^{(m-m_1)-1/2}}\label{eq1dbt}
\eea
The proof is identical as that of Lemma \ref{lemmafeasibility} except that the error estimation in (3) of Lemma \ref{best} is used. 
\end{remark}

\begin{proof}
Consider the Legendre series 
$$(x_r)^{(m_1)}(t) \sim \dsum_{n=0}^{N-r-m_1+1} a_n^{N}(m_1)L_n(t)$$
A sequence of polynomials $x_1^N(t), \cdots, x_{r+m_1}^N(t)$ is defined as follows,
\bean
x_{r+m_1}^N(t) &=& \dsum_{n=0}^{N-r-m_1+1} a_n^{N}(m_1)L_n(t)\\
x_{r+m_1-1}^N(t)&=&(x_r)^{(m_1-1)}(-1) +\displaystyle{\int_{-1}^t}x_{r+m_1}^N(s)ds\\
&\vdots&\\
x_{r+1}^N(t)&=&\dot x_r(-1) +\displaystyle{\int_{-1}^t}x_{r+2}^N(s)ds\\
\eean
and 
\bean
x_{i}^N(t)&=&x_i(-1) +\displaystyle{\int_{-1}^t}x_{i+1}^N(s)ds, \;\;\mbox{ for } 1\leq i\leq r
\eean
Define
$$x^N(t)=\MT x_1^N(t)&\cdots&x^N_r(t)\EM^T$$
and define
$$u^N(t)=\Fr{x_{r+1}^N(t)-f(x^N(t))}{g(x^N(t))}$$
From the definition of $x^N(t)$, we have $x^N(-1)=x(-1)$. If $m_1\geq 2$, then $x_{r+1}(-1)=\dot x_r(-1)$. From the definition of $u^N(t)$, we know $u^N(-1)=u(-1)$ provided $m_1\geq 2$. Therefore, $(x^N(t),u^N(t))$ satisfies (\ref{eq1b2}).
It is obvious that $x_i^N(t)$ is a polynomial of degree less than or equal to $N$; and $(x^N(t), u^N(t))$ satisfies the differential equation (\ref{eqad3}). Because we assume $V(x_r^{(m)})<\infty$, we have $x_r^{(m)}\in L^2$. From Lemma \ref{best}
\bean
||x_{r+m_1}^N(t)-x_r^{(m_1)}(t)||_\infty &=&|| x_r^{(m_1)}(t)-\dsum_{n=0}^{N-r-m_1+1}a_n^{N}(m_1)L_n(t)||_\infty \\
&\leq & C_1 || x_r||_{W^{m,2}}(N-r-m_1+1)^{-(m-m_1)+3/4}
\eean
for some constant $C_1>0$. Therefore, 
\bea
|x^N_{r+m_1-1}(t)-(x_r)^{(m_1-1)}(t)| &\leq& \displaystyle{\int_{-1}^t}|x^N_{r+m_1}(s)- (x_r)^{(m_1)}(s)|ds\nn\\
&\leq & 2C_1|| x_r||_{W^{m,2}}(N-r-m_1+1)^{-(m-m_1)+3/4}\nn
\eea
Similarly, we can prove (\ref{eq1a}) and (\ref{eq1d}). 

To prove (\ref{eq1c}), note that the spectral coefficient $a^{N}_n(m_1)$ of $(x_r^N)^{(m_1)}(t)$ is the same as the spectral coefficients of $(x_{r})^{(m_1)}(t)$. From Jackson's Theorem (Lemma \ref{Jackson}), we have
$$
|a^{N}_n(m_1)|< \Fr{6}{\sqrt{\pi}}(U(x^{(m_1+1)}_r)+V(x^{(m_1+1)}_r))\Fr{1}{n^{3/2}}
$$

In a bounded set around $x(t)$, we have $g(x)>\alpha>0$ for some $\alpha>0$ because $f$ and
$g$ are Lipschitz continuous (Definition of Problem B). Therefore, the function 
$$\Fr{s-f(x)}{g(x)}$$
is Lipschitz in a neighborhood of $(x,s)$, i.e. there exists
a constant $C_2$ independent of $N$ such that 
\bea | u^N(t) - u(t) | & = &
\left| \frac{x_{r+1}^N(t)-f(x^N(t))}{g(x^N(t)) }-
\frac{\dot x_r(t)-f(x(t))}{g(x(t))} \right|\nn\\
& \leq & C_2(|x_{r+1}^N(t) - \dot x_r(t)|+|x_1^N(t)-x_1(t)|+\cdots+|x_r^N(t)- x_r(t)|)\label{eq3}
\eea
Hence, (\ref{eq1db}) follows (\ref{eq1a}), (\ref{eq1d}) and (\ref{eq3}) when $l=0$. Similarly, we can prove (\ref{eq1b}) for $l\leq L$.
\end{proof}

Now, only after this lengthy work of preparation, we are ready to prove Theorem \ref{th_rate}.

{\it Proof of Theorem \ref{th_rate}}: Let $(x^\ast (t), u^\ast (t))$ be an optimal solution to Problem B. According to Lemma \ref{lemmafeasibility} and Remark \ref{rem1}, for any positive integer $N$ that is large enough, there exists a pair of functions $(\hat x^{N}(t),\hat u^N(t))$ in which $\hat x^{N}(t)$ consists of polynomials of degree less than or equal to $N$. Furthermore, the pair satisfies the differential equation with initial conditions in Problem B and 
\bea
||\hat x^N(t)-x^\ast (t)||_\infty &<& \Fr{L}{(N-r-m_1+1)^{m-m_1-1/2}}\label{eq3_1a}\\
||\hat u^N (t)-u^\ast(t)||_\infty &<& \Fr{L}{(N-r-m_1+1)^{m-m_1-1/2}}\label{eq3_1}\\
||(\hat x_r^{N}(t))^{(l)}-(x^{\ast}_r(t))^{(l)}||_\infty &<& \Fr{L}{(N-r-m_1+1)^{m-m_1-1/2}}, \;\;\; 1\leq l\leq m_1\label{eq3_1b}
\eea
If we define
$$\hat \bu^{Nk}=\hat u^N(t_k), \hat \bx^{Nk}=\hat x^N(t_k)$$
Then $\{(\hat \bx^{N},\hat \bu^{N})\}$ satisfies (\ref{eq3_2}) and (\ref{eq3_3}) (Lemma \ref{lemma1} and \ref{lemmafeasibility}). Because $\hat x_r^N(t)$ is a polynomial of degree less than or equal to $N$ and because of (\ref{eqDxN}), we know $(\hat x^N_r(t))^{(j)}$ equals the interpolation polynomial of $\hat \bx^N_r (D^T)^{j}$. So, 
$$
\MT 1&0& \cdots &0\EM D^j (\hat \bx_r^N)^T = (\hat x^N_r(t))^{(j)}|_{t=-1}
$$
Therefore, (\ref{eq3_1b}) implies (\ref{eq3_4a}) if the bounds $\underline{\bb}_j$ and $\bb_j$ are large enough. In addition, the spectral coefficients of $\hat \bx^N_r (D^T)^{m_1}$ is the same as the spectral coefficients of $(\hat x^N_r(t))^{(m_1)}$. From (\ref{eq1c}), (\ref{zetafun}) and (\ref{eqdbound}), we have
$$\dsum_{n=0}^{N-r-m_1+1}|a^{N}_n(m_1)| \leq \bd $$
So, the spectral coefficients of $(\hat x_r^N)^{(m_1)}$ satisfies (\ref{eq3_4}). Because we select $\underline{\bb}$ and $\bar \bb$ large enough so that the optimal trajectory of the original continuous-time problem is contained in the interior of the region, then (\ref{eq3_1a}) and (\ref{eq3_1}) imply (\ref{eq3_4b}) for $N$ large enough. In summary, we have proved that $(\hat \bx^N, \hat \bu^N)$ is a discrete feasible trajectory satisfying all constraints, (\ref{eq3_2})-(\ref{eq3_4}), in Problem \BN .

Given any bounded control input $u(\cdot)$, because the system is globally Lipschitz, it uniquely determines the trajectory $x(\cdot)$ if the initial state is fixed. Therefore, the cost $J(x^\ast(\cdot),u^\ast(\cdot))$ can be considered as a functional, denoted by ${\cal J}(u)$. Because all the functions in Problem B are $C^{m}$ with $m\geq 2$, we know that ${\cal J}(u)$ has second order Fr\'echet derivative. By Lemma \ref{lemmafrechet}
\bea
\label{eq3_5}
&&|J(x^\ast(\cdot),u^\ast(\cdot))-J(\hat x^N(\cdot),\hat u^N(\cdot))|\nn\\
&=&|{\cal J}(u^\ast)-{\cal J}(\hat u^N)|\nn\\
&\leq&C_1(||u^\ast-\hat u^N||_{W^{m_1-1,\infty}}^2) \nn\\
&\leq & \Fr{C_2}{(N-r-m_1+1)^{2m-2m_1-1}}\label{eq3_6}
\eea
for some constant numbers $C_1$ and $C_2$. The last inequality is from (\ref{eq3_1}).

Now, consider $F(\hat x^N(t), \hat u^N(t))$ as a function of $t$. Let $F^N(t)$ represent the polynomial interpolation of this function at $t=t_0,t_1,\cdots,t_N$. Let $\hat p(t)$ be the best polynomial approximation of $F(\hat x^N(t), \hat u^N(t))$ under the norm of $L^\infty [-1, 1]$. Then we have
\bea
&&|J(\hat x^N(\cdot),\hat u^N(\cdot))-\bar J^N(\hat \bx^N,\hat \bu^N)|\nn\\
&=&|J(\hat x^N(\cdot),\hat u^N(\cdot))-\dsum_{k=0}^N F(\hat \bx^{Nk},\hat \bu^{Nk})w_k-E(\hat \bx^{N0},\hat \bx^{NN})| \nn\\
&=& \left| \int_{-1}^1 F(\hat x^N(t),\hat u^N(t))dt-\int_{-1}^1 F^N(t)dt\right|\nn\\
&\leq & \int_{-1}^1 |F(\hat x^N(t),\hat u^N(t))-F^N(t)|dt\nn\\
&\leq& 2(1+\Lambda_N)||\hat p(t) - F(\hat x^N(t),\hat u^N(t))||_\infty\label{eq3_7}
\eea
where 
\EQ
\label{eq3_18}
\Lambda_N \leq \Fr{2}{\pi}log(N+1)+0.685\cdots
\EE
is the Lebesgue constant. The inequality (\ref{eq3_7}) is a corollary of Lemma \ref{lebesgue}. Because $f(\cdot)$, $g(\cdot)$, and $F(\cdot)$ are $C^{m}$, it is known (Lemma \ref{best}) that the best polynomial approximation satisfies
$$||\hat p(t)-F(\hat x^N(t),\hat u^N(t))||_\infty\leq \Fr{C_3}{N^{m_1-1}} ||F(\hat x^N(t),\hat u^N(t))||_{W^{m_1-1, \infty}}$$
Because of Lemma \ref{lemma1.5}, $\{ ||F(\hat x^N(t),\hat u^N(t))||_{W^{m_1-1, \infty}} | N\geq N_1\}$ is bounded. Therefore, 
\bea
&&|J(\hat x^N(\cdot),\hat u^N(\cdot))-\bar J^N(\hat \bx^{N},\hat \bu^{N})|\leq  \Fr{(1+\Lambda_N)C_4}{N^{m_1-1}}\leq  \Fr{C_5}{N^{\alpha}}\label{eq3_8}
\eea 
for some constant numbers $C_4$ and $C_5$ independent of $N$ and any $\alpha < m_1-1$. Let 
\EQ
\label{eq3_17}
\{ (\bar x^{\ast N},\bar u^{\ast N})\}_{N=N_0}^{\infty}
\EE
be a sequence of optimal discrete solutions. Its interpolation is denoted by $(x^{\ast N}(t), u^{\ast N}(t))$. Then, similar to the derivation above, we can prove
\bea
&&|J(x^{\ast N}(\cdot),u^{\ast N}(\cdot))-\bar J^{N}(\bar x^{\ast N},\bar u^{\ast N})|\nn\\
&\leq & 2(1+\Lambda_N)||p^{N}(t)-F(x^{\ast N}(t),u^{\ast N}(t))||_\infty\label{eq3_9}\\
&\leq & \Fr{C_6(1+\Lambda_N)}{N^{m_1-1}}||F(x^{\ast N}(t),u^{\ast N}(t))||_{W^{m_1-1,\infty}}\nn
\eea 
where $p^{N}(t)$ is the best polynomial approximation of $F(x^{\ast N}(t),u^{\ast N}(t))$ with degree less than or equal to $N$. Because of Lemma \ref{lemma1.5},  $||F(x^{\ast N}(t),u^{\ast N}(t)) ||_{W^{m_1-1,\infty}} | N\geq N_1\}$ is bounded. 
So 
\bea
&&|J(x^{\ast N}(\cdot),u^{\ast N}(\cdot))-\bar J^{N}(\bar x^{\ast N},\bar u^{\ast N})|\leq  \Fr{C_7}{N^{\alpha}}\label{eq3_11}
\eea
for some constant $C_7>0$. Now, we are ready to piece together the puzzle of inequalities and finalize the proof. 
$$\begin{array}{rcllll}
&&J(x^\ast (\cdot),u^\ast(\cdot)) \\
&\leq& J(x^{\ast N}(\cdot),u^{\ast N}(\cdot))&  \left( \begin{array}{ll}
(x^{\ast N}(t),u^{\ast N}(t)) \mbox{ is a feasible }\\
\mbox{trajectory (Lemma \ref{lemma1})}\end{array}\right)\\
&\leq& \bar J^{N}(\bar x^{\ast N}, \bar u^{\ast N})+\Fr{C_7}{N^{\alpha}} & \left( \mbox{ inequality } (\ref{eq3_11})\right)\\
&\leq & \bar J^{N}(\hat x^{N},\hat u^{N})+ \Fr{C_7}{N^{\alpha}} & \left (\begin{array}{ll} (\hat x^{N},\hat u^{N}) \mbox{ is a feasible discrete }\\ \mbox{trajectory and }(\bar x^{\ast N}, \bar u^{\ast N}) \mbox{ is optimal}\end{array}\right)\\
&\leq& J(\hat x^{N}(\cdot),\hat u^{N}(\cdot)) + \Fr{C_5}{N^{\alpha}}+\Fr{C_7}{N^{\alpha}} & \left( \mbox{ inequality } (\ref{eq3_8})\right)\\
&\leq & J(x^\ast (\cdot), u^\ast (\cdot))+\Fr{C_2}{(N-r-m_1-1)^{2m-2m_1-1}}\\ &&+\Fr{C_5}{N^{\alpha}}+\Fr{C_7}{N^{\alpha}} & \left(\mbox{ inequality } (\ref{eq3_6})\right)
\end{array}$$
Therefore, 
$$0\leq J(x^{\ast N}(\cdot),u^{\ast N}(\cdot))-J(x^\ast (\cdot),u^\ast(\cdot))\leq \Fr{C_2}{(N-r-m_1-1)^{2m-2m_1-1}}+\Fr{C_5}{N^{\alpha}}+\Fr{C_7}{N^{\alpha}}$$
This inequality implies (\ref{eq3_12b}). Furthermore, (\ref{eq3_12b}) and (\ref{eq3_11}) imply (\ref{eq3_12}). 
\hfill $\Box$

According to Theorem \ref{th_rate}, the convergence rate of the approximate cost is determined by two terms with the rates
\bea
&&\Fr{1}{(N-r-m_1-1)^{2m-2m_1-1}}\sim \Fr{1}{N^{2m-2m_1-1}}\label{eqrate1}
\eea
and 
\bea
&&\Fr{1}{N^\alpha} \sim \Fr{1}{N^{m_1-1}}\label{eqrate2}
\eea
where $m$, the smoothness of $x^\ast(t)$, is fixed. However, $m_1$ can be selected provided $f(\cdot)$, $g(\cdot)$, and $F(\cdot)$ are smooth enough. Note that increasing $m_1$ will increase the rate defined by (\ref{eqrate2}), but decrease the rate defined by (\ref{eqrate1}). There is a value of $m_1$ that determines the maximum rate. Given any real number $a \in \R$, let $[a]$ be the greatest integer less than or equal to $a$. 


\begin{corollary}
\label{corollary1}
Under the same assumption as Theorem \ref{th_rate}, the convergence rate of $J( x^{\ast N}(\cdot), u^{\ast N}(\cdot))$ and $\bar J^{N}(\bar x^{\ast N},\bar u^{\ast N})$ is
$$O\left(\frac{1}{N^{[\frac{2m}{3}]-\delta}}\right)$$
in which 
$$\delta = \left\{\begin{array}{ll}
1 &0< \gamma < \frac{2}{3}\\
3\left(1-\gamma \right)&\gamma \geq \frac{2}{3}\\
1-\mbox{any positive number}, & \gamma = 0
\end{array}\right.$$
where $\gamma=\frac{2m}{3}-[\frac{2m}{3}]$. To achieve this rate,  
$$m_1= \left\{\begin{array}{lll}
\left[ \Fr{2m}{3}\right], & 0\leq \gamma < \frac{2}{3}\\
\\
\left[ \Fr{2m}{3}\right]+1, & \gamma \geq \frac{2}{3}
\end{array}\right.
$$ 
\end{corollary}

\begin{proof}
The optimal convergence rate is determined by 
$$\max_{2\leq m_1\leq m-1}\min \{ 2m-2m_1-1,\, m_1-1\}$$
The maxmin is achieved at $$m_1=\frac{2m}{3}$$ 
However, it may not be an integer. If $m_1$ is not an integer, we have two options,
$$m_1=[\frac{2m}{3}] \mbox{ or } [\frac{2m}{3}]+1$$ 
If we define 
$$\gamma=\frac{2m}{3}-[\frac{2m}{3}]$$ 
then either $m_1=\frac{2m}{3}-\gamma$ or $m_1=\frac{2m}{3}-\gamma +1$. It is straightforward to verify that
$$\min \{2m-2m_1-1,\, m_1-1\}=\left\{\begin{array}{lll}
\frac{2m}{3}-\gamma-1, & m_1=\frac{2m}{3}-\gamma\\
\frac{2m}{3}-\gamma-3(1-\gamma), & m_1=\frac{2m}{3}-\gamma+1
\end{array}\right.$$
Therefore, $\frac{2m}{3}-1-\gamma$ is larger when $\gamma < \frac{2}{3}$, and $\frac{2m}{3}-1-2(1-\gamma)$ is larger if $\gamma \geq \frac{2}{3}$. The special case at $\gamma=0$ is because of (\ref{eqrate2}) when $\frac{2m}{3}$ equals an integer. 
\end{proof}

Different from numerical computations of differential equations, solving an optimal control problem requires the approximation, (\ref{eqad3a}), of the integration as an addition to the approximation, (\ref{eq3_2}), of the differential equation. The contributions of these approximations to the overall approximation error are different; and the errors are inversely related to each other. The following theorem indicates that the rate (\ref{eqrate1}) is due to the approximation error of the differential equation and the rate (\ref{eqrate2}) is due to the approximation error of the quadrature integration rule (\ref{JN}). To verify this fact, we define the following discretization problem with exact integration.
\vspace{0.1in}

\noindent{\bf Problem \BN (J)}\ \ 
Find $\bar x\NK\in \R^{r}$ and $\bar u\NK\in \R$, $k \ = \ 0,1,\ldots,N$, that minimize
\bea J(x^N(\cdot),u^N(\cdot))&  = & \int_{-1}^{1}F(x^N(t), u^N(t))\ dt + E(x^N(-1),x^N(1)) \label{JNJ}
\eea
subject to
\bea
&&\left\{ \begin{array}{rcl}
D ( \bar x_1^N)^T &=& (\bar x_2^N)^T \\
D (\bar x_2^N)^T&=&(\bar x_3^N)^T \\
& \vdots &\\
D (\bar x_{r-1}^N)^T&=&(\bar x_r^N)^T \\
D (\bar x_r^N)^T& = & \MT f(\bar x^{N0})+g(\bar x^{N0})\bar u^{N0}\\ \vdots \\
  f(\bar x^{NN})+g(\bar x^{NN})\bar u^{NN} \EM \\
  \end{array}\right.\label{eq3_2J}\\
  \nn\\
&&\bar x^{N0}=x_0\label{eq3_3J}\\
\nn\\
 &&\underline{\bb}  \leq  \MT \bar x\NK \\ \bar u\NK\EM
\ \leq \ \bar \bb, \;\; \;\; \mbox{ for all } 0\leq k\leq N\label{eq3_4bJ}\\
&&\underline{\bb}_j \leq  \MT 1&0& \cdots &0\EM D^j (\bx_r^N)^T \ \leq \bar \bb_j, \;\; 1\leq j\leq m_1-1 \label{eq3_4aJ}\\
&& \dsum_{n=0}^{N-r-m_1+1}|a^{N}_n(m_1) |\leq \bd \label{eq3_4J}
\eea

In Problem \BN (J), $(x^N(t), u^N(t))$ is the interpolation of $(\bar x^N, \bar u^N)$. In this discretization, we approximate the differential equation by the PS method. However, the integration in the cost function is exact. In this case, the overall error is controlled by the single rate (\ref{eqrate1}) rather than the two-rate convergence of Problem \BN . Without the integration error of the cost function, the convergence rate is improved to $\Fr{1}{N^{2m-3}}$; and the smoothness requirement can be reduced to $m\geq 2$.  

\begin{theorem}
\label{th_rate2}
Suppose Problem B has an optimal solution $(x^\ast (t), u^\ast (t))$ in which the strong derivative $(x_r^{\ast}(t))^{(m)}$ has bounded variation for some $m\geq 2$. In Problem \BN (J), select $m_1$ so that $1\leq m_1\leq m-1$. Suppose $f(\cdot)$, $g(\cdot)$, and $F(\cdot)$ are $C^{m}$. Suppose all other bounds in Problem \BN\ are large enough.  
Given any sequence
\EQ
\label{eq3_13j}
\{ (\bar x^{\ast N},\bar u^{\ast N})\}_{N\geq N_1}
\EE
of optimal solutions of Problem \BN (J). Then the cost of (\ref{eq3_13j}) converges to the optimal cost at the following rate
\bea
&&\left| J(x^\ast(\cdot),u^\ast (\cdot))-J( x^{\ast N}(\cdot), u^{\ast N}(\cdot))\right|
\leq \Fr{M_1}{(N-r-m_1-1)^{2m-2m_1-1}}\label{eq3_12bj}\\
\nn
\eea
for some constants $M_1$ independent of $N$.
\end{theorem}

\begin{proof}
Let 
\EQ
\label{eq3_17a}
\{ (\bar x^{\ast N},\bar u^{\ast N})\}_{N=N_0}^{\infty}
\EE
be a sequence of optimal solutions of Problem \BN (J). According to Lemma \ref{lemmafeasibility} and Remark \ref{rem1}, for any positive integer $N$ that is large enough, there exists a pair of functions $(\hat x^{N}(t),\hat u^N(t))$ in which $\hat x^{N}(t)$ consists of polynomials of degree less than or equal to $N$. Furthermore, the pair satisfies the differential equation in Problem B and the inequalities (\ref{eq3_1a}), (\ref{eq3_1}), and (\ref{eq3_1b}). If we define
$$\hat \bu^{Nk}=\hat u^N(t_k), \hat \bx^{Nk}=\hat x^N(t_k)$$
Then, from the first part in the proof of Theorem \ref{th_rate},  $(\hat \bx^N, \hat \bu^N)$ is
a discrete feasible solution satisfying all constraints in Problem \BN (J). By Lemma \ref{lemmafrechet}
\bea
\label{eq3_5a}
&&|J(x^\ast(\cdot),u^\ast(\cdot))-J(\hat x^N(\cdot),\hat u^N(\cdot))|\nn\\
&=&|{\cal J}(u^\ast)-{\cal J}(\hat u^N)|\nn\\
&\leq&C_1(||u^\ast-\hat u^N||_{W^{m_1-1,\infty}}^2) \nn\\
&\leq & \Fr{C_2}{(N-r-m_1+1)^{2m-2m_1-1}}\label{eq3_6a}
\eea
for some constant numbers $C_1$ and $C_2$. The last inequality is from (\ref{eq1db}). The interpolation $(x^{\ast N}(t), u^{\ast N}(t))$ of (\ref{eq3_17a}) is a feasible trajectory of Problem B (Lemma \ref{lemma1}). Thus,
$$\begin{array}{rcllll}
&&J(x^\ast (\cdot),u^\ast(\cdot)) \\
&\leq& J(x^{\ast N}(\cdot),u^{\ast N}(\cdot))&  \left( \begin{array}{ll}
(x^{\ast N}(t),u^{\ast N}(t)) \mbox{ is a feasible }\\
\mbox{trajectory (Lemma \ref{lemma1})}\end{array}\right)\\
&\leq & J(\hat x^{N}(\cdot),\hat u^{N}(\cdot)) & \left (\begin{array}{ll} (\hat \bx^{N},\hat \bu^{N}) \mbox{ is a feasible discrete }\\ \mbox{trajectory and }(\bar x^{\ast N}, \bar u^{\ast N}) \mbox{ is optimal}\end{array}\right)\\
&\leq & J(x^\ast (\cdot), u^\ast (\cdot))+\Fr{C_2}{(N-r-m_1-1)^{2m-2m_1-1}}  & \left(\mbox{ inequality } (\ref{eq3_6a})\right)
\end{array}$$
Therefore,
$$0\leq J(x^{\ast N}(\cdot),u^{\ast N}(\cdot))-J(x^\ast (\cdot),u^\ast(\cdot))\leq \Fr{C_2}{(N-r-m_1-1)^{2m-2m_1-1}}$$
\end{proof}

\section{Existence and Convergence of Approximate Optimal Solutions}
\label{sec-convergence}
\setcounter{equation}{0}
In Section \ref{rate}, the rate of convergence for the cost function is proved. However, the results do not guarantee the convergence of the approximate optimal trajectory $\{(x^N(t), u^N(t))\}$. In this section, we prove the existence of feasible trajectories for Problem \BN\ and the existence of a convergent subsequence in any set of approximate optimal solutions. In addition, we consider a larger family of problems. Different from Section \ref{formulation} where Problem B does not contain constraints other than the control system, in this section the problem of optimal control may contain nonlinear path constraints. Furthermore, general endpoint conditions are allowed, rather than being limited to the initial value problem as in the previous sections. 
\vspace{0.1in}

\noindent{\bf Problem B:}\ \ 
Determine the state-control function pair $(x(t),u(t))$, $x\in \R^r$ and $u\in \R $, that minimizes the cost function
\bea
J(x(\cdot),u(\cdot))&  = & \int_{-1}^{1}F(x(t), u(t))\ dt + E(x(-1),x(1)) \label{Joriginal}
\eea
subject to the state equation
\bea
&&\left\{\begin{array}{lll}
\dot x_1 = x_2\\
\;\;\; \vdots  \\
\dot x_{r-1} =  x_r\\
\dot x_r = f(x) + g(x)u 
\end{array}\right.
\label{state-eq}
\eea
end-point conditions
\bea
e(x(-1),x(1)) & = & 0 \label{ini}
\eea
and state-control constraints
\bea
h(x(t),u(t) ) & \leq  & 0 \label{cons-c}
\eea
where $x\in \R^r$, $u\in \R $, and $F: \R^r \times \R \to \R$, $E: \R^r \times
\R^r \to \R$, $f: \R^r \to \R$, $g: \R^r \to \R$ $e: \R^r \times \R^r \to
\R^{N_e}$ and $h: \R^r \times \R \to \R^{N_h}$ are all Lipschitz continuous functions with respect to their arguments. In addition, we assume $g(x)\neq 0$ for all $x$. The corresponding discretization is defined as follows.
\vspace{0.1in}

\noindent{\bf Problem ${\bf B}^{\bf N}$:}\ \ 
Find $\bar x\NK\in \R^{r}$ and $\bar u\NK\in \R$, $k \ = \ 0,1,\ldots,N$, that minimize
\bea \bar J^N(\bar x^N ,\bar u^N) & = & \sum_{k=0}^{N}F(\bar x\NK,\bar u\NK)w_k + 
E(\bar x^{N0},\bar x^{NN}) \label{JN}
\eea
subject to
\bea
&&\left\{ \begin{array}{rcl}
D ( \bar x_1^N)^T &=& (\bar x_2^N)^T \\
D (\bar x_2^N)^T&=&(\bar x_3^N)^T \\
& \vdots &\\
D (\bar x_{r-1}^N)^T&=&(\bar x_r^N)^T \\
D (\bar x_r^N)^T& = & \MT f(\bar x^{N0})+g(\bar x^{N0})\bar u^{N0}\\ \vdots \\
  f(\bar x^{NN})+g(\bar x^{NN})\bar u^{NN} \EM \\
  \end{array}\right.\label{state-N}\\
  \nn\\
&&\|e(\bar x^{N0},\bar x^{NN})\|_\infty  \leq  (N-r-1)^{-\beta} \label{ini-N}\\
&&h(\bar x\NK,\bar u\NK)  \leq  (N-r-1)^{-\beta}\cdot \mathbf{1}, \qquad \ \ \ \ 
 \mbox{ for all } 0\leq k\leq N\label{bounds}\\
 \nn\\
 &&\underline{\bb}  \leq  \MT \bar x\NK \\ \bar u\NK\EM
\ \leq \ \bar \bb, \;\; \;\; \mbox{ for all } 0\leq k\leq N\label{bound2}\\
&&\underline{\bb}_j \leq  \MT 1&0& \cdots &0\EM D^j (\bx_r^N)^T \ \leq \bar \bb_j, \mbox{ if } 1\leq j\leq m_1-1 \mbox{ and } m_1\geq 2\label{bound3b}\\
&& \dsum_{n=0}^{N-r-m_1+1}|a^{N}_n(m_1) |\leq \bd \label{bound3}
\eea

The discretization is almost identical to the one used in the previous sections except for the path constraints and the endpoint conditions, which must be treated with care. Note that in Problem \BN the right sides of (\ref{ini-N}) and ({\ref{bounds}) are not zero. It is necessary to relax (\ref{ini}) and (\ref{cons-c}) by a small margin for the reason of feasibility. The margin approaches zero as $N$ is increased. Without this relaxation, it is shown by a counter example in \cite{gong} that Problem \BN\ may have no feasible trajectories. 

Some feasibility and convergence results were proved in \cite{gong}, which take the form of consistent approximation theory based on the convergence assumption about $\{\dot x_r^N(t)\}$ and $\{\bx^{N0}\}$. The goal of this section is to remove this bothersome assumption by using a fundamentally different approach. In addition, the proofs in this section are not based on necessary conditions of optimal control and any coercivity assumption, which are widely used in existing work on the convergence of direct optimal control methods. Before we introduce main results in this section, some useful results from \cite{gong} are summarized in the following Lemma.
\begin{lemma} (\cite{gong}) \label{T-2}
Suppose Problem B has an optimal solution $(x^\ast (t), u^\ast (t))$ satisfying $x^\ast_r(t)\in W^{m,\infty}$, $m\geq 2$. Let $\{(\bar x^N,\bar u^N)\}_{N=N_1}^\infty$ be a sequence of feasible solutions to \refBN. Suppose there is a subsequence $\{N_j\}_{j=1}^\infty$ of $\{ N\}_{N=1}^{\infty}$ such that the sequence $\left\{\bar x^{N_j0}\right\}_{j=1}^{\infty}$ converges as $N_j\rightarrow\infty$. Suppose there exists a continuous function $q(t)$ such that $\dot x_r^{N_j}(t)$ converges to $q(t)$ uniformly in $[-1,1]$. Then, there exists $(x^\infty(t), u^\infty(t))$ satisfying \refB\, such that the following limits converge uniformly in $[-1, 1]$.
\bea
&&\lim_{N_j\rightarrow \infty}(x^{N_j}(t) - x^\infty(t)) =  0 \label{con3}\\
&&\lim_{N_j\rightarrow \infty}(u^{N_j}(t) -u^\infty(t)) =  0 \label{con4} \\
&&\lim_{N_j\rightarrow\infty} \bar J^{N_j} (\bx^{N_j},\bu^{N_j})  = 
J(x(\cdot),u(\cdot))\label{converge-j}\\
&&\lim_{N_j\rightarrow\infty} J (x^{N_j},u^{N_j})  = 
J(x(\cdot),u(\cdot))\label{converge-jb}
\eea
In addition to the above assumptions, if $\{(\bar x^N,\bar u^N) \}_{N=N_1}^\infty$ is a sequence of optimal solutions subject to the constraints \refBN, then $(x^\infty(t), u^\infty(t))$ must be an optimal solution to Problem B.
\end{lemma}

The following are the two main theorems of this section. Relative to \cite{kang}, these results has a tightened bounds for $m_1$ and $\beta$.

\begin{theorem}\ \ \label{maintheorem}
(Existence of solutions) Consider Problem B and Problem \BN\ defined in Section \ref{formulation}. Suppose Problem B has a feasible trajectory $(x (t),u (t))$ in which $(x_r (t))^{(m)}$ has bounded variation for some $m\geq 2$. In Problem \BN , let $m_1$ be any integer and $\beta$ be any real number satisfying $1\leq m_1\leq m-1$ and $0<\beta < (m-m_1)-\frac{3}{4}$. Then, there exists $N_1>0$ so that, for all $N\geq N_1$, Problem \BN\ has a feasible trajectory satisfying (\ref{state-N})-(\ref{bound3}). Furthermore, $(x (t),u (t))$ and the interpolation $(x^N(t),u^N(t))$ satisfy (\ref{eq1a})-(\ref{eq1c}).
\end{theorem} 

\begin{theorem}
\label{maintheorem2}
(Convergence) Consider Problem B and Problem \BN\ defined in Section \ref{formulation}.
Suppose Problem B has an optimal solution $(x^\ast (t),u^\ast (t))$ in which $(x^\ast_r (t))^{(m)}$ has bounded variation for some $m\geq 3$. In Problem \BN , let $m_1$ be any integer and $\beta$ be any real number satisfying $2\leq m_1\leq m-1$ and $0<\beta < (m-m_1)-\frac{3}{4}$. Then for any sequence $\{(\bar x^{\ast N},\bar u^{\ast N})\}_{N=N_1}^\infty$ of optimal solutions of Problem \BN , there exists a subsequence, $\{(\bar x^{\ast N_j},\bar u^{\ast N_j})\}_{j\geq 1}^\infty$, and an optimal solution, $(x^{\ast}(t), u^{\ast}(t))$, of Problem B so that the following limits converge uniformly in $[-1, 1]$
\bea
\lim_{N_j\rightarrow\infty} (x^{ \ast N_j}(t)- x^\ast (t)) & = & 0\nn\\
\lim_{N_j\rightarrow\infty} (u^{ \ast N_j}(t)- u^\ast(t) )& = & 0\label{eqtheoremii}\\
\lim_{N_j\rightarrow\infty} \bar J^{N_j}(\bx^{\ast N_j}, \bu^{\ast N_j})&=& J(x^\ast(\cdot),u^\ast(\cdot))\nn\\
\lim_{N_j\rightarrow\infty}  J(x^{\ast N_j}(\cdot), u^{\ast N_j}(\cdot))&=& J(x^\ast(\cdot),u^\ast(\cdot))\nn
\eea 
where $(x^{\ast N_j}(t),u^{\ast N_j}(t))$ is the interpolation of $(\bar x^{\ast N},\bar u^{\ast N})$.   
\end{theorem}

\begin{remark}
The integers $m$ and $m_1$ in Theorem \ref{maintheorem} are smaller than those in Theorem \ref{maintheorem2}, i.e. the existence theorem is proved under a weaker smoothness assumption than the convergence theorem. 
\end{remark}

To prove these theorems, we first briefly review some results on real analysis and then prove a lemma.  Given a sequence of functions $\{ f_k(t)\}_{k=1}^{\infty}$ defined on $[a,b]$. It is said to be uniformly equicontinuous if for every $\epsilon > 0$, there exists a $\delta > 0$ such that for all $t$, $t'$ in $[a,b]$ with $|t' - t| < \delta$, we have 
$$|f_k(t) - f_k(t')| < \epsilon$$ 
for all $k\geq 1$. The following Proposition and Theorem are standard in real analysis \cite{sansone}. 

\begin{proposition}
\label{propappendix}
If $f_k(t)$ is differentiable for all $k$, and if $\{ \dot f_k(t)\}_{k=1}^{\infty}$ is bounded. Then, $\{ f_k(t)\}_{k=1}^{\infty}$ is uniformly equicontinuous.
\end{proposition}

\begin{theorem} (Arzel\`{a}-Ascoli Theorem) Consider a sequence of continuous functions $\{ h_n(t)\}_{n=1}^{\infty}$ defined on a closed interval $[a,b]$ of the real line with real values. If this sequence is uniformly bounded and uniformly equicontinuous, then it admits a subsequence  which converges uniformly. 
\end{theorem}

\begin{lemma}
\label{lemma2}
(\cite{kang})Let $\{(\bar x^N,\bar u^N)\}_{N=N_1}^\infty$ be a sequence satisfying \refBN. Assume the set 
\bea
&&\left\{ \left. ||\ddot x^N_r(t)||_\infty \right| N\geq N_1\right\}\label{eqlm3}
\eea
is bounded. Then, there exists $(x^\infty(t), u^\infty(t))$ satisfying \refB\ and a subsequence $\{(\bx^{N_j},\bu^{N_j})\}_{N_j\geq N_1}^{\infty}$ such that (\ref{con3}), (\ref{con4}), (\ref{converge-j}) and (\ref{converge-jb}) hold. Furthermore, if $\{(\bar x^N,\bar u^N) \}_{N=N_1}^\infty$ is a sequence of optimal solutions to Problem \BN, then $(x^\infty(t), u^\infty(t))$ must be an optimal solution to Problem B.
\end{lemma}

\begin{proof} 
Let $x_r^N(t)$ be the interpolation polynomial of $\bar x_r^N$. Because (\ref{eqlm3}) is a bounded set, we know that the sequence of functions $\{ \dot x_r^N(t) | N\geq N_1\}$ is uniformly equicontinuous (Proposition \ref{propappendix}). By the Arzel\`{a}-Ascoli Theorem, a subsequence $\{ \dot x_r^{N_j}(t)\}$ converges uniformly to a continuous function $q(t)$. In addition, because of (\ref{bound2}), we can select the subsequence so that $\{ \bx^{N_j0}\}_{N_j\geq N_1}^{\infty}$ is convergent. Therefore, all conclusions in Lemma \ref{T-2} hold true. 
\end{proof}

Now, we are ready to prove the theorems.
\vspace{0.1in}

{\it Proof of Theorem \ref{maintheorem}:} For the feasible trajectory $(x(t),u(t))$, consider the pair $(x^N(t),u^N(t))$ in Lemma \ref{lemmafeasibility} that satisfies the differential equation (\ref{state-eq}). Define 
\bea
\label{eq5}
\begin{array}{rcl}
\bar x\NK&=&x^N(t_k)\\
\bar u\NK&=&u^N(t_k)
\end{array}
\eea 
for $0\leq k\leq N$. From Lemma \ref{lemma1}, we know that $\{(\bar x^N, \bar u^N)\}$ satisfies the discrete equations in (\ref{state-N}). In the next we prove that the mixed state-control constraint (\ref{bounds}) is
satisfied. Because $h$ is Lipschitz continuous and because of (\ref{eq1a}) and (\ref{eq1db}), there exists a constant $C$ independent of $N$ so that
\bea \| h(x(t),u(t)) - h( x^N(t), u^N(t)) \|& \leq &
C(|x_1(t)-x^N_1(t)|
+\cdots+|x_r(t)-x^N_r(t)|+ | u(t) -u^N(t) | )\nn\\
& \leq & CMV|| x_r||_{W^{m,2}}(r+1)(N-r-m_1+1)^{-(m-m_1)+3/4}\nn
\eea
Hence
\bea h(x^N(t),u^N(t)) & \leq & h(x(t),u(t)) +
CM|| x_r||_{W^{m,2}}(r+1)(N-r-m_1+1)^{-(m-m_1)+3/4}
\cdot\mathbf{1}\nn\\
& \leq & CM|| x_r||_{W^{m,2}}(r+1)(N-r-m_1+1)^{-(m-m_1)+3/4} \nn \eea 
Because $\beta < m-m_1-\frac{3}{4}$, there exists a positive integer $N_1$ such that, for all $N>N_1$, 
\bea
CM|| x_r||_{W^{m,2}}(r+1)(N-r-m_1+1)^{-(m-m_1)+3/4} & \leq & (N-r-1)^{-\beta} \nn \eea
Therefore $x^N_1(t_k)$, $\ldots$, $x^N_r(t_k)$, $u^N(t_k)$,
$k=0,1,\ldots, N$, satisfy the mixed state and control constraint
(\ref{bounds}) for all $N>N_1$. 

By a similar procedure, we can prove that the endpoint condition (\ref{ini-N}) is satisfied. Because $x_r^N(t)$ is a polynomial of degree less than or equal to $N$, and because of (\ref{eqDxN}) and (\ref{eq5}), we know $(x^N_r(t))^{(j)}$ equals the interpolation polynomial of $\bx^N_r (D^T)^{j}$. So, 
$$
\MT 1&0& \cdots &0\EM D^j (\bx_r^N)^T = (x^N_r(t))^{(j)}|_{t=-1}
$$
Therefore, (\ref{eq1d}) implies (\ref{bound3b}) if the interval between $\underline{\bb}_j$ and $\bb_j$ is large enough. In addition, the spectral coefficients of $\bx^N_r (D^T)^{m_1}$ is the same as the spectral coefficients of $(x^N_r(t))^{(m_1)}$. From (\ref{eq1c}) and (\ref{eqdbound}), we have
$$\dsum_{n=0}^{N-r-m_1+1}|a^{N}_n(m_1)| \leq \bd $$
So, $\{(\bar x^N,\bar u^N)\}$ satisfies (\ref{bound3}). Because we select $\underline{\bb}$ and $\bar \bb$ large enough so that the optimal trajectory of the original continuous-time problem is contained in the interior of the region, we can assume that $(x(t),u(t))$ is also bounded by $\underline{\bb}$ and $\bar \bb$. Then, (\ref{eq1a}) and (\ref{eq1db}) imply (\ref{bound2}) for $N$ large enough. To summarize, $(\bar x^N, \bar u^N)$ satisfies (\ref{state-N})-(\ref{bound3}). Therefore, it is
a feasible trajectory of \refBN.
\hfill $\Box$
\vspace{0.1in}

{\it Proof of Theorem \ref{maintheorem2}:}
Consider $\{(\bar x^{\ast N},\bar u^{\ast N})\}_{N=N_1}^\infty$, a sequence of optimal solutions of Problem \BN. From Lemma \ref{lemma1.5}, 
$$\{ ||\ddot x_r^{\ast N}(t)||_\infty | N\geq N_1\}$$
is bounded. Now, we can apply Lemma \ref{lemma2} to conclude that there exists a subsequence of $\{(x^{\ast N}(t),u^{\ast N}(t))\}_{N=N_1}^\infty$ and an optimal solution of Problem B so that the limits in (\ref{eqtheoremii}) converge uniformly.
\hfill $\Box$

\section{Simulation results}
The rate of convergence for the optimal cost is illustrated in the following example
\bean
&&\min_u \int_0^\pi (1-x_1 +x_1x_2+x_1u)^2dt\\
&&\mbox{subject to}\\
&&\dot x_1=-x_1^2x_2\\
&&\dot x_2=-1+\frac{1}{x_1}+x_2+\sin t +u\\
&&x(0)=\MT 1\\0\EM,\;\; x(\pi)=\MT \frac{1}{\pi+1}\\ 2\EM
\eean
The analytic solution of this problem is known so that the approximation error can be computed 
\bean
&&x_1(t)=\frac{1}{1-\sin t +t}\\
&&x_2(t)=1-\cos t\\
&&u=-(t+1)+\sin t +\cos t\\
&&\mbox{optimal cost}=0
\eean
The problem is solved using PS optimal control method. The approximated optimal cost is compared to the true value. The number of nodes, N, ranges from $4$ to $16$. The error decreases rapidly as shown in Table \ref{table1}.
\begin{table}[h]
	\centering
		\begin{tabular}{|c|c|c|c|c|c|c|c|}
			\hline
			N&4&6&8&10&12&14&16\\
			\hline
Error &$7.5\times10^{-2}$&$1.1\times 10^{-3}$&$2.1\times 10^{-4}$&$7.1\times 10^{-5}$&$6.7\times 10^{-6}$&$1.0\times 10^{-6}$&$5.8\times 10^{-7}$\\
			\hline
		\end{tabular}
		\label{table1}
		\caption{The error of optimal cost}
		\label{table2}
\end{table}
The rate of convergence is illustrated in the following Figure \ref{figad}.
\begin{figure}[h]
	\centering
		\includegraphics[width=3in]{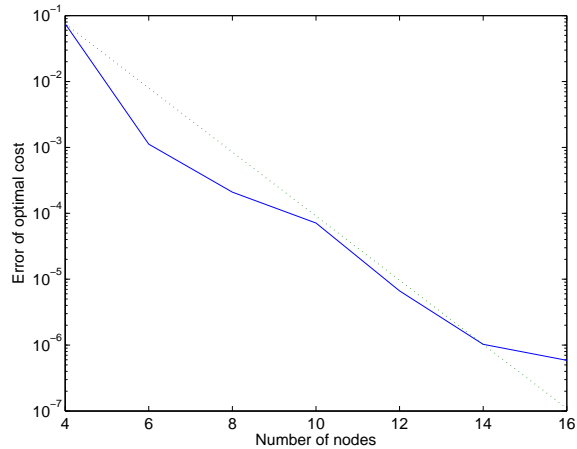}
		\caption{Log-scale plot of the error of optimal cost (the solid curve)}
		\label{figad}
\end{figure}
Because the analytic solution is $C^\infty$, the rate of convergence of the PS method is faster than any polynomial rate. As a result, it converges exponentially. Of course, in practical computations the accuracy is limited by the machine precision of the computers. Therefore, the accuracy cannot be improved after $N$ is sufficiently large.

\section{Conclusions}
\setcounter{equation}{0}
It is proved that the PS optimal control has a high-order rate of convergence. According to the theorems in Section \ref{rate}, the approximate cost computed using the Legendre PS method converges at an order determined by the smoothness of the original problem. More specifically, the rate is about $\frac{1}{N^{2m/3-1}}$, where $m$ is defined by the smoothness of the optimal trajectory. If the cost function can be accurately computed, then the convergence rate is improved to $\frac{1}{N^{2m-1}}$. If the optimal control is $C^\infty$, then the convergence rate can be made faster than any given polynomial rate. The results in Section \ref{sec-convergence} imply that the discretization using the Legendre PS method is feasible; and there always exists a convergent subsequence from the approximate discrete optimal solutions, provided some smoothness assumptions are satisfied.


\begin{thebibliography}{99}

\bibitem{betts:book}
J. T. Betts, {\it Practical Methods for Optimal Control Using Nonlinear
Programming}, SIAM, Philadelphia, PA, 2001.
\bibitem{betts:survey}
J. T. Betts,   ``Survey of Numerical Methods for Trajectory Optimization,''
 {\it Journal of Guidance, Control, and Dynamics,} Vol. 21, No. 2,
1998, pp. 193-207.
\bibitem{boyd} J. P. Boyd, Chebyshev and Fourier Spectral Methods, second edition, Dover, 2001,
\bibitem{canuto} C. Canuto, M. Y. Hussaini, A. Quarteroni and T. A. Zang, {\em Spectral Method in Fluid Dynamics.} New York: Springer-Verlag, 1988.
\bibitem{elnagar2}
G. Elnagar, and  M. A. Kazemi, Pseudospectral Chebyshev Optimal Control of
Constrained Nonlinear Dynamical Systems, {\it Computational Optimization and
Applications,} 11, 1998, pp. 195-217.
\bibitem{fahroo} Fahroo, F., Ross, I. M., "Costate Estimation by a Legendre Pseudospectral Method," Proceedings of the AIAA Guidance, Navigation and Control Conference, 10-12 August 1998, Boston, MA.
\bibitem{gong} Q. Gong, W. Kang, and I. M. Ross, A Pseudospectral Method for the Optimal Control of Constrained Feedback Linearizable Systems, IEEE Trans. Automat. Contr., Vol. 51, No. 7, pp. 1115-1129, 2006.
\bibitem{gong2} Q. Gong, M. Ross, W. Kang, F. Fahroo, Connections Between the Covector Mapping Theorem and Convergence of Pseudospectral Methods for Optimal Control, Computational Optimization and Applications, to appear.
\bibitem{H2000} W. W. Hager, Runge-Kutta methods in optimal control and the transformed adjoint system, {\em Numerische Mathematik,} Vol. 87, pp. 247-282, 2000.
\bibitem{DH2000}  A. L. Dontchev and W. W. Hager, The Euler approximation in state constrained
optimal control, {\em Mathematics of Computation,} Vol. 70, pp. 173-203, 2000.
\bibitem{hesthaven} J. Hesthaven, S. Gottlieb, and D. Gottlieb, Spectral Methods for Time-Dependent Problems, Cambridge University Press, 2007.
\bibitem{kangbedrossian} W. Kang, N. Bedrossian, Pseudospectral Optimal Control Theory Makes Debut Flight - Saves NASA \$1M in under 3 hrs, SIAM News, September, 2007. 
\bibitem{kang2} W. Kang, Q. Gong, and I. M. Ross, On the Convergence of Nonlinear Optimal Control using Pseudospectral Methods for Feedback Linearizable Systems, International Journal of Robust and Nonlinear Control, Vol. 17, 1251-1277, online publication, 3 January, 2007.
\bibitem{kang} W. Kang, I. M. Ross, Q. Gong, Pseudospectral Optimal Control and Its Convergence Theorems, Analysis and Design of Nonlinear Control Systems - In Honor of Alberto Isidori, A. Astolfi and L. Marconi eds., Springer, 2008.  
\bibitem{otis} S. W. Paris and C. R. Hargraves, \textit{OTIS 3.0 Manual}, Boeing Space and
Defense Group, Seattle, WA, 1996.
\bibitem{paris:aas06} S. W. Paris, J. P. Riehl, and W. K. Sjauw, ``Enhanced Procedures for
Direct Trajectory Optimization Using Nonlinear Programming and
Implicit Integration,'' {\it Proceedings of the AIAA/AAS
Astrodynamics Specialist Conference and Exhibit}, 21-24 August 2006,
Keystone, CO.  AIAA Paper No. 2006-6309.
\bibitem{polak:book}
E. Polak, \textit{Optimization: Algorithms and Consistent Approximations},
Springer-Verlag, Heidelberg, 1997.
\bibitem{riehl:aas06} J. P. Riehl, S. W. Paris, and W. K. Sjauw, ``Comparision of Implicit
Integration Methods for Solving Aerospace Trajectory Optimization
Problems,'' {\it Proceedings of the AIAA/AAS Astrodynamics
Specialist Conference and Exhibit}, 21-24 August 2006, Keystone, CO.
AIAA Paper No. 2006-6033.
\bibitem{ross} Ross, I. M., A Beginner's Guide to DIDO: A MATLAB Application Package for Solving
Optimal Control Problems, Elissar Inc., Monterey, CA, October 2007.
\bibitem{sansone} G. Sansone, A. H. Diamond, and E. Hille, Orthogonal Functions, Robert E. Krieger Publishing Co., Huntington, New York, 1977. 
\bibitem{wouk} A. Wouk, A Course of Applied Functional Analysis, John Wiley \& Sons, New York, 1979.

\end{thebibliography}
\end{document}